
\newdimen\firstident
\newdimen\secondident
\newdimen\thirdident
\newbox\identbox
\newcount\assertioncount


\font\fractur=cmfrak




\let\pz=\S

%

\let\eps=\varepsilon

%

\def\Abf{{\bf A}}
\def\Sbf{{\bf S}}

%

\def\Bbar{{\bar B}}

\def\hbar{{\bar h}}

\def\Mbar{{\bar M}}

\def\sbar{{\bar s}}

\def\xbar{{\bar x}}

\def\ybar{{\bar y}}

%

%

%

%

\def\Bline{\underline{B}}

\def\hline{\underline{h}}

\def\Mline{\underline{M}}

\def\Nline{\underline{N}}

\def\Sline{\underline{S}}

\def\Xline{\underline{X}}

%

\def\ftilde{{\tilde f}}

\def\Gtilde{{\tilde G}}

\def\mtilde{{\tilde m}}
\def\Mtilde{{\tilde M}}

\def\Ntilde{{\tilde N}}

\def\utilde{{\tilde u}}
\def\Utilde{{\tilde U}}

\def\Xtilde{{\tilde X}}

%

\def\Blinebar{{\overline{\Bline}}}

\def\Mlinebar{{\overline{\Mline}}}

\def\Slinebar{{\overline{\Sline}}}

%

\def\Xlinetilde{{\tilde{\Xline}}}

%


%

%

\def\Qdbar{{\bar {\Q}}}

%

\def\mgd{\mu\mkern-9mu\mu}

%

\def\Gdhat{{\hat {\G}}}

%

%

\def\hgbar{{\bar \eta}}

%

%

%

%

%

\def\Oscrhat{{\hat \Oscr}}

%

%

\def\Mtttilde{{\widetilde \Mtt}}

%

%

\def\Sscrdot{{\dot \Sscr}}


\catcode`\º=\active\defº{\int}
\catcode`\¶=\active\def¶{\partial}
\catcode`\Æ=\active\defÆ{\triangle}
\catcode`\Â=\active\defÂ{{\neg}}
\catcode`\Ã=\active\defÃ{{\sqrt}}
\catcode`\Å=\active\defÅ{{\approx}}
\catcode`\°=\active\def°{{\infty}}
\catcode`\­=\active\def­{\mathrel{\not=}}
\catcode`\²=\active\def²{\leq}
\catcode`\³=\active\def³{\geq}
\catcode`\Ö=\active\defÖ{\div}
\catcode`\É=\active\defÉ{{\ldots}}
\catcode`\Ç=\active\defÇ{\ll}
\catcode`\¯=\active\def¯{\emptyset}
\catcode`\±=\active\def±{\pm}
\catcode`\È=\active\defÈ{\gg}
\catcode`\¡=\active\def¡{\circ}
\def\star{{}^*}

\def\vdual{{}^{\mathord{\vee}}}               

\def\dlbrack{\lbrack{\mskip-2mu}\lbrack}      
\def\drbrack{\rbrack{\mskip-2mu}\rbrack}      
\def\lrangle{\langle\>{,}\>\rangle}           

\def\Qp{\Q_{p}}

\def\Qpbar{{{\overline {\Q}}_{p}}}

\def\Zp{\Z_{p}}


\def\overneq#1#2{\lower0.5pt\vbox{\lineskiplimit\maxdimen\lineskip-.5pt
                \ialign{$#1\hfil##\hfil$\crcr#2\crcr\not=\crcr}}}


\def\stimes{\mathbin{\raise1pt\hbox{$\scriptscriptstyle \bf
            \vert$}\mkern-5mu\times}}   


\def\ad{\mathop{\rm ad}\nolimits}

\def\Aut{\mathord{\rm Aut}}
\def\Autline{\mathop{\underline {\Aut}}\nolimits}
\def\bi{{\rm bi}}

\def\codim{\mathop{\rm codim}\nolimits}

\def\Def{\mathop{\rm Def}\nolimits}

\def\det{\mathop{\rm det}\nolimits}

\def\dim{\mathop{\rm dim}\nolimits}
\def\End{\mathop{\rm End}\nolimits}

\def\et{{\rm \acute et}}

\def\Gal{\mathop{\rm Gal}\nolimits}

\def\Hom{\mathop{\rm Hom}\nolimits}
\def\Homline{\mathop{\underline {\Hom}}\nolimits}

\def\id{\mathop{\rm id}\nolimits}
\def\Im{\mathop{\rm Im}\nolimits}

\def\inf{\mathop{\rm inf}}

\def\int{\mathop{\rm int}\nolimits}

\def\inv{\mathop{\rm inv}\nolimits}
\def\Isom{\mathop{\rm Isom}\nolimits}
\def\Isomline{\mathop{\underline {\Isom}}\nolimits}

\def\Ker{\mathop{\rm Ker}\nolimits}

\def\Lie{\mathop{\rm Lie}\nolimits}
\def\log{\mathop{\rm log}\nolimits}

\def\mult{\mathop{\rm mult}\nolimits}

\def\Res{\mathop{\rm Res}\nolimits}

\def\Spec{\mathop{\rm Spec}\nolimits}
\def\Spf{\mathop{\rm Spf}\nolimits}

\def\Tr{\mathop{\rm Tr}\nolimits}

\def\setback(#1){\mathrel{\mkern-#1mu}}
\def\varfill#1{$\smash{#1} \mkern-8mu
   \cleaders\hbox{$\mkern-3mu \smash{#1} \mkern-3mu$}\hfill
   \mkern-8mu #1$}
\def\equalfill{\varfill{=}}


\let\ar=\rightarrow

\let\air=\hookrightarrow

\let\asr=\mapsto

\let\arr=\longrightarrow
\let\all=\longleftarrow
\def\airr{\lhook\joinrel\arr}          

\def\arriso{\buildrel\sim\over\arr}    
\def\arrover#1{\buildrel#1\over\arr}   

\def\airrover#1{\buildrel#1\over\airr} 




\def\allover#1{\buildrel#1\over\all}   

\def\arvar(#1){\hbox to #1pt{\rightarrowfill}}
\def\alvar(#1){\hbox to #1pt{\leftarrowfill}}
\def\arvarover(#1)#2{\mathop{\arvar(#1)}\limits^{#2}}
\def\alvarover(#1)#2{\mathop{\alvar(#1)}\limits^{#2}}


\let\ad=\downarrow
\let\au=\uparrow
\def\add{\Big\ad}                      
\def\auu{\Big\au}                      
\def\addleft#1{\llap{$\vcenter{\hbox{$\scriptstyle #1$}}$}\add}
\def\auuleft#1{\llap{$\vcenter{\hbox{$\scriptstyle #1$}}$}\auu}
\def\addright#1{\add\rlap{$\vcenter{\hbox{$\scriptstyle #1$}}$}}
\def\auuright#1{\auu\rlap{$\vcenter{\hbox{$\scriptstyle #1$}}$}}






\def\aqrvar(#1){\hbox to #1pt{\equalfill}}







%
%

\def\set#1#2{\{\,#1\>\vert\>\hbox{#2}\,\}}

%
%

%
%

\def\powerseries over #1 in #2{{#1 \dlbrack #2 \drbrack}}

%
%

\def\laurentseries over #1 in #2{{#1 (\!( #2 )\!)}}

%
%

\def\smallmatrix(#1,#2;#3,#4){\left({{#1\atop #3}\>{#2\atop #4}}\right)}

%
%

%
%

\def\tensor#1{\otimes_{#1}}

%
%

%
%

%

%

\def\restricted#1{\vert_{#1}}

%
%

\def\subscript#1\atop#2{_{\scriptstyle #1 \atop #2}}

%
%

%
%




\newfam\frac
\textfont\frac=\fractur

\def\mfr{{\fam\frac m}}

%
%
%
%
\font\tenss=cmss10
\newfam\ssfam %
\textfont\ssfam=\tenss
\catcode`\_=11
\def\suf_fix{}
\def\scaled_rm_box#1{%
 \relax
 \ifmmode
   \mathchoice
    {\hbox{\tenrm #1}}%
    {\hbox{\tenrm #1}}%
    {\hbox{\sevenrm #1}}%
    {\hbox{\fiverm #1}}%
 \else
  \hbox{\tenrm #1}%
 \fi}
\def\suf_fix_def#1#2{\expandafter\def\csname#1\suf_fix\endcsname{#2}}
\def\I_Buchstabe#1#2#3{%
 \suf_fix_def{#1}{\scaled_rm_box{I\hskip-0.#2#3em #1}}
}
\def\rule_Buchstabe#1#2#3#4{%
 \suf_fix_def{#1}{%
  \scaled_rm_box{%
   \hbox{%
    #1%
    \hskip-0.#2em%
    \lower-0.#3ex\hbox{\vrule height1.#4ex width0.07em }%
   }%
   \hskip0.50em%
  }%
 }%
}
\I_Buchstabe B22
\rule_Buchstabe C51{34}
\I_Buchstabe D22
\I_Buchstabe E22
\I_Buchstabe F22
\rule_Buchstabe G{525}{081}4
\I_Buchstabe H22
\I_Buchstabe I20
\I_Buchstabe K22
\I_Buchstabe L20
\I_Buchstabe M{20em }{I\hskip-0.35}
\I_Buchstabe N{20em }{I\hskip-0.35}
\rule_Buchstabe O{525}{095}{45}
\I_Buchstabe P20
\rule_Buchstabe Q{525}{097}{47}
\I_Buchstabe R21 
\rule_Buchstabe U{45}{02}{54}
\suf_fix_def{Z}{\scaled_rm_box{Z\hskip-0.38em Z}}
\catcode`\"=12
\newcount\math_char_code
\def\suf_fix_math_chars_def#1{%
 \ifcat#1A
  \expandafter\math_char_code\expandafter=\suf_fix_fam
  \multiply\math_char_code by 256
  \advance\math_char_code by `#1
  \expandafter\mathchardef\csname#1\suf_fix\endcsname=\math_char_code
  \let\next=\suf_fix_math_chars_def
 \else
  \let\next=\relax
 \fi
 \next}
%
%
%
%
\def\font_fam_suf_fix#1#2 #3 {%
 \def\suf_fix{#2}
 \def\suf_fix_fam{#1}
 \suf_fix_math_chars_def #3.
}
\font_fam_suf_fix
 0rm
 ABCDEFGHIJKLMNOPQRSTUVWXYZabcdefghijklmnopqrstuvwxyz
\font_fam_suf_fix
 2scr
 ABCDEFGHIJKLMNOPQRSTUVWXYZ
\font_fam_suf_fix
 \slfam sl
 ABCDEFGHIJKLMNOPQRSTUVWXYZabcdefghijklmnopqrstuvwxyz
\font_fam_suf_fix
 \bffam bf
 ABCDEFGHIJKLMNOPQRSTUVWXYZabcdefghijklmnopqrstuvwxyz
\font_fam_suf_fix
 \ttfam tt
 ABCDEFGHIJKLMNOPQRSTUVWXYZabcdefghijklmnopqrstuvwxyz
\font_fam_suf_fix
 \ssfam
 ss
 ABCDEFGHIJKLMNOPQRSTUVWXYZabcdefgijklmnopqrstuwxyz
\catcode`\_=8
%
%
%
%
%
%
%
%
%
%
%
%
%
%
%
%
%
%
%
%
%
%
%
%
%
%



\def\en_item#1#2{%
 \par
 \setbox\identbox=\hbox #1{#2}%
 \noindent
 \hangafter=1%
 \hangindent=\wd\identbox
 \box\identbox
 \ignorespaces
}

\def\ennopar_item#1#2{%
 \setbox\identbox=\hbox #1{#2}%
 \noindent
 \hangafter=1%
 \hangindent=\wd\identbox
 \box\identbox
 \ignorespaces
}

%
%

\def\hangitem#1{\en_item{}{#1\enspace}}

%
%
%

%

\def\hanghangitem#1{\en_item{}{\kern\firstident #1\enspace}}

%
%
\def\hanghanghangitem#1{\en_item{}{\kern\secondident #1\enspace}}


\def\newpage{\vfil\eject}

%
%
\def\indention#1{%
 \setbox\identbox =\hbox{\kern\parindent{#1}\enspace}%
 \firstident=\wd\identbox
}

%
%
\def\subindention#1{%
 \setbox\identbox=\hbox{{#1}\enspace}%
 \secondident=\wd\identbox
 \advance \secondident by \firstident
}

%
%
\def\subsubindention#1{%
 \setbox\identbox=\hbox{{#1\ }\enspace}%
 \thirdident=\wd\identbox
 \advance \thirdident by \secondident
}%

%
%
\def\litem#1{\en_item{to\firstident}{\kern\parindent#1\hfil \enspace }}
\def\llitem#1{\en_item{to\secondident}{\kern\firstident#1\hfil\enspace}}
\def\lllitem#1{\en_item{to\thirdident}{\kern\secondident#1\hfil\enspace}}
\def\ritem#1{\en_item{to\firstident}{\kern\parindent\hfil#1\enspace}}
\def\rritem#1{\en_item{to\secondident}{\kern\firstident\hfil#1\enspace}}
\def\rrritem#1{\en_item{to\thirdident}{\kern\secondident\hfil#1\enspace}}
\def\citem#1{\en_item{to\firstident}{\kern\parindent\hfil#1\hfil\enspace}}
\def\ccitem#1{\en_item{to\secondident}{\kern\firstident\hfil#1\hfil\enspace}}
\def\cccitem#1{\en_item{to\thirdident}{\kern\secondident\hfil#1\hfil\enspace}}
\def\rmlitem#1{\en_item{to\firstident}{\kern\parindent{\rm #1}\hfil \enspace }}
\def\rmllitem#1{\en_item{to\secondident}{\kern\firstident{\rm #1}\hfil\enspace}}
\def\rmlllitem#1{\en_item{to\thirdident}{\kern\secondident{\rm #1}\hfil\enspace}}

%
%

\def\lnoitem#1{\ennopar_item{to\firstident}{\kern\parindent#1\hfil \enspace }}
\def\llnoitem#1{\ennopar_item{to\secondident}{\kern\firstident#1\hfil\enspace}}
\def\lllnoitem#1{\ennopar_item{to\thirdident}{\kern\secondident#1\hfil\enspace}}

%
%

\def\assertionlist{
   \assertioncount=1
   \indention{(2)}}

\def\assertionitem{
   \litem{{\rm(\number\assertioncount)}}\advance\assertioncount by 1
}

%
%

\def\bulletlist{
   \indention{$\bullet$}}

\def\bulletitem{
   \litem{$\bullet$}}

%
%


\newcount\parno
\newcount\secno
\newcount\subsecno
\newdimen\partisize
\newdimen\parindentvalue


\magnification 1200

\hsize=156 true mm
\vsize=8.9 true in

\tolerance=300
\pretolerance=100

\parindentvalue=15pt
\parindent=\parindentvalue

\mathsurround=0pt

\normallineskiplimit=.5pt
\normalbaselineskip=15pt
\normallineskip=1pt plus .5 pt minus .5 pt
\normalbaselines

\abovedisplayskip = 8pt plus 3pt minus 3pt
\abovedisplayshortskip = 1pt plus 2pt
\belowdisplayskip = 8pt plus 3pt minus 3pt
\belowdisplayshortskip = 5pt plus 2pt

\partisize=12pt


\global\parno=0
\global\secno=0
\global\subsecno=0


\newskip\parskipamount               
 \parskipamount=21pt plus4pt minus4pt
 \def\parskip{\vskip\parskipamount}
\newskip\secskipamount               
 \secskipamount=18pt plus4pt minus2pt
 \def\secskip{\vskip\secskipamount}


\font\titlefont=cmbx12 at 18pt
\font\partifont=cmbx7 at \partisize




\outer\def\paragraph#1{\goodbreak
   \global\secno=0
   \global\advance\parno by 1
   \parskip
   \hangitem{\partifont \the\parno~}
   {\partifont #1}
   \par\nobreak}

%

\outer\def\Notation#1: #2\par{\nobreak\secskip
  \itemitem{}{\bf Notation#1}:\enspace{\sl#2\par}\penalty50}

%

\outer\def\secstart#1{\penalty-100
   \global\subsecno=0
   \global\advance\secno by 1
   \secskip\noindent
   {\bf (\the\parno.\the\secno)~#1}}

%

\outer\def\proof#1: {\smallbreak {\it Proof\ #1}\/:\enspace}

%

\outer\def\claim#1: #2\par{\smallbreak
   {\bf #1}:\enspace{\sl#2\par}\penalty50}


\outer\def\bibliography{
   \parskip
   \parindent=0pt
   {\partifont Bibliography}
   \parskip
   \frenchspacing
   }



%

\def\lformno{
   \global\advance\subsecno by 1
   \leqno(\the\parno.\the\secno.\the\subsecno)
   }


\def\displayno{
   \global\advance\subsecno by 1
   (\the\parno.\the\secno.\the\subsecno)
   }


\def\actualsecno{(\the\parno.\the\secno)}


\def\actualsubsecno{(\the\parno.\the\secno.\the\subsecno)}

%
%
%

\def\label#1{\xdef#1{\actualsecno}}

%
%
%

\def\sublabel#1{\xdef#1{\actualsubsecno}}

\def\pisogi{p\mathord{\hbox{-}}{\tt Isog}_{\Iscr,K^p}}
\def\pisog{p\mathord{\hbox{-}}{\tt Isog}}


\centerline{\titlefont Congruence relations for Shimura varieties}
\secskip
\centerline{\titlefont associated to some unitary groups}
\parskip
\centerline{Oliver B\"ultel}
\centerline{Torsten Wedhorn}

\bigskip\bigskip\bigskip


{\partifont Introduction}

\secskip

Let $G/\Q$ be a reductive group and $X$ be a $G(\R)$-conjugacy class of
homomorphisms $\Res_{\C/\R}(\G_m) \arr G_{\R}$ such that $(G,X)$ is a
Shimura datum and denote for some compact open subgroup $K \subset
G(\Abf_f)$ by ${\rm Sh}_K(G,X)$ the associated Shimura variety. Its
$\C$-valued points are given by $G(\Q)\backslash ((X \times G(\Abf_f)/K)$
and it is defined over a number field $E$ which is associated to
$(G,X)$. Hence its cohomology carries an action of $\Gamma_E$ (the absolute
Galois group of $E$) and of the Hecke algebra $\Hscr(G(\Abf_f)//K)$ of
$K$-biinvariant locally constant $\Q$-valued functions with compact
support, and these actions commute with each other. It is one of the
fundamental problems to understand the precise relation of these two
actions. One step towards this understanding is the following: Let $p$ be a
prime of good reduction (i.e.\ $G_{\Q_p}$ is unramified and we have $K =
K_pK^p$ with $K_p \subset G(\Q_p)$ hyperspecial) and let $v$ be a place
of $E$ lying over $p$. Denote by ${\rm Fr}_v$ the conjugacy class of
corresponding geometric Frobenius elements in $\Gamma_E$. Blasius and
Rogawski [BR] have defined a polynomial $H = H_{(G,X),p}$ (the so-called
Hecke polynomial) with coefficients in the local
Hecke algebra $\Hscr(G(\Q_p)//K_p)$ and they conjectured that $H({\rm
Fr}_v) = 0$. 

This is known for Shimura curves for a long time by the work
of Eichler/Shimura and Ihara. Further, Faltings and Chai [FC] have proved
in the Siegel case that ${\rm Fr}_v$ satisfies an equation of the correct
degree. Further the first of us has proved the
conjecture in the case of certain orthogonal groups [??], and the second in the
PEL-case when $G_{\Qp}$ is a split reductive group [We2].

The goal of this article is to prove this conjecture in the case where $G$
is the group of $B$-linear similitudes of a skew-hermitian space
$(V,\lrangle)$ where $B$ is a simple $\Q$-algebra endowed with a positive
involution which is split at the chosen prime $p$ which we assume to be
odd. We further suppose that $G_{\R}$ is isomorphic to a group auf unitary
similitudes of a hermitian form of signature $(n-1,1)$. We finally assume
that $n$ is even (in particular, the case of Picard modular surfaces is
excluded). In fact, we also assume that $B_{\Qp}$ is simple, as otherwise
we are in the situation of [We2]. See chapter 1 for the precise
assumptions.

As we are in the PEL-case, ${\rm Sh}_K(G,X)$ can be written as a moduli
space of abelian varieties with certain additional structures
and it has a model $\Mscr$ over $O_{E_{(p)}}$. To prove the congruence
relation we proceed in three steps: First we study the reduction of $\Mscr$
to characteristic $p$. In the second part we examine a moduli space of
isogenies which is used to define the Hecke correspondence and an extension
in characteristic $p$. In the final chapter we use the obtained results to
prove the congruence relation.

We are now going to describe the results of the three parts in more
detail. For this we assume only that the signature is of the form $(n-1,1)$,
but $n$ can be arbitrary (the assumption ``$n$ even'' is needed only in the
last chapter).

In the first part (chapter 3 -- 5) we look at three stratifications of the
special fibre $\Mscr_0$ of $\Mscr$. These stratifications are defined as
follows: Geometric points of $\Mscr_0$ are given by isomorphism classes of
tuples $(A,\iota,\lambda,\hgbar)$ where $A$ is an abelian variety, $\iota$
is an action of an order $O_B$ of $B$, $\lambda$ is a
$p$-principal polarization, and $\hgbar$ is a level structure. The strata
are given by
\indention{(b)}
\litem{(a)} The isomorphism class of the isocrystal associated to
$(A,\iota,\lambda)$ (the Newton polygon stratification),
\litem{(b)} The isomorphism class of the $p$-torsion of $(A,\iota,\lambda)$
(the Ekedahl-Oort stratification),
\litem{(c)} The isomorphism class of the $p$-divisible group of
$(A,\iota,\lambda)$ (the final stratification).

Our main result in this part is the following

\smallskip

{\bf Theorem}: {\sl On the non-supersingular locus of $\Mscr_0$ the three
stratifications coincide.}

\smallskip

This implies that every Newton polygon stratum is a union of Ekedahl-Oort
strata and as a corollary we get from the corresponding fact for the Ekedahl
Oort stratification [We3] that the closure of a Newton polygon stratum is
again a union of Newton polygon strata and that all Newton polygon strata
are equi-dimensional, nonempty and of the expected dimension. In particular
we see that the supersingular locus $\Mscr^{\rm ss}$ has dimension strictly
smaller than half of the dimension of $\Mscr_0$ if and only if $n$ is even
(if $n$ is odd, we have $\dim(\Mscr^{\rm ss}) = \dim(\Mscr_0)/2$).

In the second part (chapter 6 and 7 and proposition 8.7) we look at the special
fibre $\pisog_0$ of the moduli space of $p$-isogenies $f\colon
(A,\iota,\lambda,\hgbar) \ar (A',\iota',\lambda',\hgbar')$. Here our main
results are
\indention{--}
\litem{--} The canonical morphisms source and target $\pisog_0 \arr
\Mscr_0$ are finite over the non-supersingular locus and they are finite
locally free over the generic stratum (with respect to any of the three
stratifications).
\litem{--} We describe $\pisog_0$ generically using an analogue of the
local model which we call ``generic model''. In particular we find that
every component of $\pisog_0$ has dimension greater or equal than
$\dim(\Mscr_0)$.
\litem{--} Let $\pisog^{\mu}$ be the inverse image of
the generic stratum of $\Mscr_0$ in $\pisog_0$. Then every geometric point
of $\pisog^{\mu}$ admits a canonical lift to characteristic zero.

As a corollary we get that $\pisog^{\mu}$ is dense in $\pisog_0$ if
and only if $n$ is even (if $n$ is odd one can show that there
supersingular components in $\pisog_0$ and those components are even in the
closure of the generic fibre of $\pisog$).

In the third part (chapter 8) we prove the principal result. Here we assume
that $n$ is even. The main idea
of the proof is the same as in [FC] and [We2]. We consider the $\Q$-vector
space $\Q[\pisog \otimes E]$ (resp.\ $\Q[\pisog_0]$) generated by the
components of $\pisog \otimes E$ (resp.\ $\pisog_0$). We have
$\Q[\pisog^{\mu}] = \Q[\pisog_0]$ because $\pisog^{\mu}$ is dense in
$\pisog_0$ . Both $\Q$-vector spaces obtain a $\Q$-algebra structure by
composition of $p$-isogenies. We define a commutative diagram of
$\Q$-algebra homomorphisms
$$\matrix{\Hscr(G(\Qp)//K_p)_0 & \arrover{h^0} & \Q[\pisog \otimes E] \cr
\addleft{\Sscrdot} && \addright{{\rm red}_p} \cr
\Hscr(L(\Qp)//(K_p \cap L(\Qp)))_0 & \arrover{\hbar} &
\Q[\pisog_0].\cr}\leqno{(*)}$$
Here $L \subset G_{\Qp}$ is the centralizer of the norm of the minuscule
coweight $\mu$ of $G$ associated
to $(G,X)$. The subscript $(\ )_0$ in $\Hscr(G(\Qp)//K_p)_0$ denotes the
subalgebra of the Hecke algebra of functions with integral support such
that if we invert the characteristic function of $K_p(p\cdot\id_V)K_p$ we get
$\Hscr(G(\Qp)//K_p)$, similar for $L$. We take $\Q$ as the
coefficient ring for the Hecke algebras. The homomorphism $\Sscrdot$ is a
twisted version of the Satake homomorphism already defined in [We2]. It is
injective. On the right hand side ${\rm red}_p$ denotes the specialization
map of cycles. Then elements in $\Q[\pisog \otimes E]$ give rise to
correspondences on ${\rm Sh}_K(G,X)$ and $h^0$ is simply the definition of
Hecke correspondences on ${\rm Sh}_K(G,X)$. Finally, to construct $\hbar$
we define it as a homomorphism $\Hscr(L(\Qp)//(K_p \cap
L(\Qp)) \arr \Q[\pisog^{\mu}]$ and then use $\Q[\pisog^{\mu}] =
\Q[\pisog_0]$.

Once we have the diagram $(*)$, the result follows from a group theoretic
statement of elements in the Hecke algebra which has been proved by the
first author in his thesis. 

\medskip

We will now give an overview on the structure of this article: In the first
chapter we introduce notations and the second
chapter introduces the Dieudonn\'e modules with the additional structure
which are considered in the sequel, the so-called unitary Dieudonn\'e
modules of signature $(n-1,1)$.

In the third chapter we give a precise description of the isomorphism
classes of the reduction mod $p$ of the unitary Dieudonn\'e modules of
signature $(n-1,1)$ using a theorem of Moonen [Mo] which simplifies our
previous arguments considerably (3.5). Further we calculate the dimension
of the group of automorphisms of such a reduction (3.13) because this gives the
codimension of the corresponding Ekedahl-Oort stratum by [We3].

The fourth chapter contains a criterion to recover the Newton polygon of a
unitary Dieudonn\'e module from its reduction mod $p$ (4.2). In particular,
we see that the Oort stratification is finer than the Newton polygon
stratification. Further in (4.9) we give the description of the
homomorphism induced by $p$-isogenies induced on Dieudonn\'e modules for
non-supersingular points.

In the fifth chapter we compare Newton polygon, Ekedahl-Oort and final
stratification in detail and obtain the dimension of the supersingular
locus of ${\rm Sh}_K(G,X) \otimes \kappa(E_v)$.

The content of the sixth chapter is the definition of the generic model
for the moduli space of $p$-isogenies and the lower bound for the dimension
of components of $\pisog \otimes \kappa(E_v)$ (6.15).

In the seventh chapter we examine the source and target morphisms $\pisog
\otimes \kappa(E_v) \arr {\rm Sh}_K(G,X) \otimes \kappa(E_v)$. We show that
they are finite (7.3) over the non-supersingular locus. In (7.5) we
describe the deformation functors of $\mu$-ordinary $p$-isogenies which
allows to show that the $\mu$-ordinary locus in $\pisog \otimes
\kappa(E_v)$ is dense if $n$ is even (7.8). It further implies that in this
case every point of $\pisog \otimes \kappa(E_v)$ can be lifted to
characteristic zero (7.9) and that source and
target are finite locally free over the generic stratum although this last
fact is proved only in 8.7 due to reasons of notation.

The last chapter contains now the definition of $h^0$ and $\hbar$,
statement and proof of the main theorem (8.10 -- 8.12).

\medskip

We are grateful to U.\ G\"ortz and S.\ Orlik for helpful remarks on
preliminary versions of the manuscript.

\newpage

\paragraph{Notations}

\secstart{}\label{\genno} Let $p > 2$ be a fixed prime. Let $\Qdbar_p$ be
an algebraic closure of $\Q_p$ and fix an embedding $\Qdbar \air \Qdbar_p$
where $\Qdbar$ denotes the algebraic closure of $\Q$ in $\C$. Further, we fix a
square root $\sqrt{-1}$ of $-1$ in $\C$.

\secstart{}\label{\Shimurano} Throughout we
denote by $\Iscr = (B,\star,V,\lrangle,G,O_B,\Lambda,h,\mu)$ the
collection of the following data:
\bulletlist
\bulletitem $B$ denotes a simple $\Q$-algebra such that $B \tensor{\Q} \R
\cong M_m(\C)$ and such that $B \tensor{\Q} \Q_p \cong M_m(K)$ where $K$ is
a quadratic unramified extension of $\Q_p$.
\bulletitem $\star$ is a positive involution on $B$.
\bulletitem $V \not= (0)$ is a finitely generated left $B$-module.
\bulletitem $\lrangle$ is a nondegenerate skew-hermitian $\Q$-valued form
on $V$ (i.e.\ $\lrangle$ is alternating and $\langle bv,v'\rangle = \langle
v,b\star v'\rangle$ for all $b \in B$ and $v,v' \in V$).
\bulletitem $G$ denotes the algebraic $\Q$-group of $B$-linear symplectic
similitudes of $(V,\lrangle)$.
\bulletitem $O_B$ is a $\star$-invariant $\Z_{(p)}$-order such that $O_B
\otimes \Z_p$ is a maximal order of $B_{\Qp}$.
\bulletitem $\Lambda$ denotes an $O_B$-invariant $\Z_p$-lattice of $V
\tensor{\Q} \Q_p$ such that the alternating form on $\Lambda$ which is
induced from $\lrangle$ is a perfect $\Z_p$-form.
\bulletitem $h$ is a homomorphism of real algebraic groups $h\colon
\Sbf = \Res_{\C/\R}(\G_{m,\C}) \arr G_{\R}$ such that $h$ defines a Hodge
structure of type $\{(-1,0),(0,-1)\}$ on $V$ (with the sign convention of
[De]) and such that
$$V_{\R} \times V_{\R} \arr \R, \qquad (v,v') \asr \langle
v,h(\sqrt{-1})v'\rangle$$
is symmetric and positive definite on $V_{\R}$.
\bulletitem $\mu$ denotes the composition
$$\G_{m,\C} \arr \Sbf_{\C} \arrover{h_{\C}} G_{\C}$$
where the first arrow is the embedding $\G_{m,\C} \air \Sbf_{\C} =
\prod_{\Gal(\C/\R)} \G_{m,\C}$ whose image is the factor of $\Sbf_{\C}$
corresponding to the identity in $\Gal(\C/\R)$.

We can and we do fix an isomorphism $B \otimes \Q_p \cong M_m(K)$ such that
$O_B \otimes \Z_p$ is identified with $M_m(O_K)$.

\secstart{}\label{\moreno} The algebraic group $G$ over $\Q$ is reductive
and connected. The assumption $B \otimes \R \cong M_m(\C)$ implies that
$G_{\R}$ is isomorphic to the group of unitary similitudes $GU(r,s)$ of an
hermitian form of signature $(r,s)$ for nonnegative integers $r$ and $s$ with
$$r+s = n := \dim_{\Q}(V)/2m.$$
Without loss of generality we can assume that $r \geq s$.

We denote by $(X^*,R,X_*,R\vdual,\Delta)$ the based root datum of $G$ with its
action of $\Gamma = \Gal(\Qdbar/\Q)$ and by $\Omega$ its Weyl group. Let
$\{\mu\}$ be the $G(\C)$-conjugacy class of $\mu$ which we consider as an
element of $X_*/\Omega$. The group $\Gamma$ acts continuously on the
discrete groups $X_*$ and
$\Omega$ such that $\gamma(w\lambda) = {}^\gamma w\gamma(\lambda)$ for
$\gamma \in \Gamma$, $w \in \Omega$ and $\lambda \in X_*$. Hence, we get an
induced continuous action of $\Gamma$ on $X_*/\Omega$. The stabilizer of
$\{\mu\}$ in $\Gamma$ is open and defines a finite extension $E$ of $\Q$
which is called the {\it reflex field of $\Iscr$}. It is easy to see (e.g.\
[We1] 2.3.2) that $E = \Q$ if $r = s$ and that $E$ is a quadratic imaginary
extension otherwise.

Via the chosen embedding $\Qdbar \air \Qdbar_p$ we can
consider $\{\mu\}$ as a $G(\Qdbar_p)$-conjugacy class of cocharacters. Denote
by $v\vert p$ the place of $E$ given by the chosen embedding $\Qdbar \air
\Qpbar$. Then the field of definition of $\{\mu\}$ with respect to the
$\Gal(\Qdbar_p/\Qp)$-action is the $v$-adic completion of $E$. Let
$\kappa$ be the residue class field of $v$. The action of
$\Gal(\Qdbar_p/\Qp)$ factors over $\Gal(K/\Qp)$. We denote by $\sigma$ the
nontrivial element in $\Gal(K/\Qp)$.

Note that the assumption on $B_{\Qp}$ and the existence of $\Lambda$ imply
that $G_{\Q_p}$ is unramified over $\Qp$. More precisely, the stabilizer
of $\Lambda$ defines a reductive model $\Gtilde$ of $G_{\Qp}$ over $\Zp$.

\secstart{}\label{\modspaceno} Denote by $\Abf^p_f$ the ring of finite
adeles of $\Q$ with trivial $p$-th component. We fix an open compact
subgroup $K^p \subset G(\Abf^p_f)$ and denote by $\Mscr =
\Mscr_{\Iscr,K^p}$ the associated moduli space, defined by Kottwitz
[Ko]. We assume $K^p$ to be sufficiently small, such that the moduli
problem is representable by a scheme. This scheme is then smooth and
quasi-projective over the localization of $O_E$ in $p$. It classifies tuples
$(A,\iota,\lambda,\hgbar)$ where $A$ is an abelian scheme up to
prime-to-$p$-isogeny, $\lambda$ is a $\Q$-homogeneous polarization of $A$
of degree prime to $p$,
$\iota\colon O_B \ar \End(A) \otimes_{\Z} \Z_{(p)}$ is an involution
preserving $\Z_{(p)}$\nobreak-algebra homomorphism where the involution is
$\star$ on $O_B$ and the Rosati-Involution given by $\lambda$ on $\End(A)
\otimes_{\Z} \Z_{(p)}$, and where $\hgbar$ is a $K^p$-level
structure. Further $(A, \iota, \lambda, \hgbar)$ satisfies the determinant
condition, i.e.\ we have an identity of polynomial functions on $O_B$
$$\det(b \vert \Lie(A)) = \det(b \vert V_0)$$
where $V_0 \subset V_{E_v}$ is the weight zero space of some $\mu \in
\{\mu\}$ (see [Ko2] \pz5 or [RZ] 3.23 a) for a precise formulation of the
determinant condition).

It follows from the general theory of Shimura varieties that
$\Mscr_{\Iscr,K^p}$ has relative dimension
$$2\langle \rho, \{\mu\} \rangle = rs$$
where $\rho$ denotes the halfsum of positive roots in the based root datum
of $G$.

\secstart{}\label{\defineLevi} We choose in the conjugacy class $\{\mu\}$ the
unique element $\mu$ which is dominant. As $G_{E_v}$ is unramified over
$E_v$, it is already defined over $O_{E_v}$. Set
$$N(\mu) = \sum_{\gamma \in \Gal(E_v/\Qp)} \gamma(\mu).$$
This a one-parameter subgroup of $G_{\Qp}$ which extends (uniquely) to a
homomorphism $\G_{m,\Zp} \arr \Gtilde$. We denote by $L \subset G_{\Qp}$
the centralizer of $N(\mu)$ in $G_{\Qp}$. This is a Levi subgroup of
$G_{\Qp}$ which is unramified as reductive group over $\Qp$..

\secstart{}\label{\defineHecke} Denote by $K_p$ the stabilizer of
$\Lambda$ in $G(\Q_p)$ with the notations of \Shimurano, i.\ e.\ $K_p =
\Gtilde (\Z_p)$. Let $\Hscr_{\Q}(G)
= \Hscr_{\Q}(G(\Qp)//K_p)$ be the Hecke algebra of $\Q$-valued locally
constant $K_p$-biinvariant functions on $G(\Qp)$. It is commutative as
$K_p$ is a hyperspecial subgroup of $G(\Qp)$. The Cartan decomposition of
$G(\Qp)$ shows that we have an identification
$$K_p\backslash G(\Qp)/K_p =
(X_*/\Omega)^{\sigma}\lformno\sublabel{\Cartanident}$$
where $\sigma \in \Gal(K/\Qp)$ denotes the nontrivial element.

We denote by $\Hscr_{\Q}(G)_0$ the subalgebra of $\Hscr_{\Q}(G)$ which
consists of functions whose support lies in $\End_{\Z_p}(\Lambda)$. If $P$
is the double coset of $p\id_V$, we have $\Hscr_{\Q}(G)_0[1/P] =
\Hscr_{\Q}(G)$ as $p\id_V$ is central in $G(\Qp)$.

We denote by $H = H_{G,\{\mu\}} \in \Hscr_{\Q}(G)_0$ the Hecke polynomial
associated to $G$ and the conjugacy class $\{\mu\}$ (see e.g.\ [We2] \pz 2
for the definition and the fact that its coefficients lie in
$\Hscr_{\Q}(G)_0$).

Further ${}_LK = L(\Qp) \cap K_p$ is a hyperspecial subgroup of
$M(\Qp)$. It is also equal to the centralizer of the extension of $N(\mu)$
in $\Gtilde(\Zp)$. Again we get Hecke algebras $\Hscr_{\Q}(L) =
\Hscr_{\Q}(L(\Qp)//{}_LK)$ and $\Hscr_{\Q}(L)_0$.

\secstart{}\label{\definepisog} Let $S$ be an $O_{E,v}$-scheme and let
$(A_1,\iota_1,\lambda_1,\hgbar_1), (A_2,\iota_2, \lambda_2, \hgbar_2) \in
\Mscr(S)$ be two $S$-valued points. A {\it $p$-isogeny}
$$f\colon (A_1,\iota_1,\lambda_1,\hgbar_1) \arr (A_2,\iota_2,
\lambda_2,\hgbar_2)$$
is an $O_B$-linear isogeny $f\colon A_1 \arr A_2$ such that there exists an
integer $c \geq 0$ with $(\alpha f)^*\lambda_2 = p^c\lambda_1$.

The integer $c = c(f)$ is called the {\it multiplicator of $f$}.

We have the obvious notion of the pullback of a $p$-isogeny with respect to
a morphism of schemes $S' \arr S$ and the notion of an isomorphism of
two $p$-isogenies. We get a functor $\pisog = \pisogi$ on the category
of $O_{E,\nu}$-schemes. The functor $\pisog^c$
of isomorphism classes of $p$-isogenies with multiplicator $c$ is an open
and closed subfunctor.

An important example of a $p$-isogeny is the relative Frobenius
$${\rm Fr}\colon A \arr A^{(q)}\lformno\sublabel{\defineFrob}$$
where $q$ denotes the number of elements in $\kappa(E_v)$, i.e.\ $q = p$ if
$r = s$ and $q = p^2$ otherwise. Its multiplicator is $\log_pq$.

Sending a $p$-isogeny to its source (resp.\ its target) defines
morphisms
$$s\colon \pisog^c \arr \Mscr_{\Iscr,K^p}, \qquad t\colon \pisog^c \arr
\Mscr_{\Iscr,K^p}$$
which are representable by projective surjective morphisms. In particular
$\pisog^c$ (and hence $\pisog$) is representable by a scheme. Further, the
morphisms are finite \'etale in characteristic zero. Hence $\pisog \otimes
E$ is smooth.

For every $p$-isogeny $f$ we denote by $f\vdual$ its dual. It is again a
$p$-isogeny of the same multiplicator as $f$.


\paragraph{Unitary p-divisible groups}

\secstart{}\label{\defpdivmodule} Let $S$ be an $O_{E_v}$-scheme which is
$p$-adically complete. A {\it $p$-divisible $\Iscr$-module
over $S$} is a triple $(X,\iota,\lambda)$ where $X$ is a $p$-divisible group
over $S$, $\iota\colon O_K \arr \End(X)$ is a ring homomorphism and where
$\lambda\colon (X,\iota) \arriso (X\vdual,\iota\vdual)$ is an $O_K$-linear
isomorphism of $p$-divisible groups such that $\lambda\vdual = -\lambda$.
Here $(\ )\vdual$ denotes the Serre dual and by $\iota\vdual\colon O_K \arr
\End(X\vdual)$ we mean the homomorphism $\iota\vdual(a) =
\iota(a\star)\vdual$ for $a \in O_K$. We further require that we have for
all $n \geq 0$ an identity of polynomial functions on $O_K$
$$\det(\iota(a)\vert \Lie(X \times_S S_n)) = \det(a\vert V_0
\tensor{M_m(K)} K^m)\lformno\sublabel{\pdivdetcond}$$
where $S_n = \Spec(\Oscr_S/p^{n+1}\Oscr_S)$. As $X \times_S S_n$ is formally
smooth over $S_n$, $\Lie(X \times_S S_n)$ is a finite locally free
$\Oscr_{S_n}$-module and \pdivdetcond\ holds for all $n \geq 0$ if it holds
for one $n \geq 0$.

We want to reinterprete this condition further. We fix the isomorphism
$$O_K \tensor{\Zp} O_K \arriso O_K \times O_K, \qquad a \otimes a' \asr
(aa',a\star\,a').\lformno\sublabel{\unitdecomp}$$
As the determinant condition \pdivdetcond\ can be checked after some finite
\'etale base change $S' \arr S$ we can assume that $S$ is also an
$O_K$-scheme. We obtain a decomposition of the $O_K \otimes \Oscr_S$-module
$\Lie(X)$ of the form $\Lie(X) = L_0 \oplus L_1$ where $L_i$ is an
$\Oscr_S$-submodule. The determinant condition is then equivalent to the
condition
$${\rm rk}_{\Oscr_S}(L_0) = r, \qquad {\rm rk}_{\Oscr_S}(L_1) = s$$
where $r$ and $s$ are the integers defined in \moreno.

\secstart{}\label{\assocdm} Let $S$ be as above and let
$(A,\iota,\lambda,\hgbar)$ be an $S$-valued point of $\Mscr_{\Iscr,K^p}$. The
$p$-primary torsion of $(A,\iota,\lambda)$ defines a $p$-divisible
$O_B \otimes \Z_p$-module $(X',\iota')$ together with an $O_B$-linear
isomorphism $\lambda'\colon (X',\iota') \arriso
(X'\vdual,\iota'\vdual)$. As $O_B \otimes \Zp = M_m(O_K)$ we can apply
Morita equivalence and we get a $p$-divisible $\Iscr$-module over $S$ which we
call {\it associated to $(A,\iota,\lambda,\hgbar)$}.

\secstart{} We are going to study the category of $p$-divisible
$\Iscr$-modules over a fixed algebraically closed field $k$ of characteristic
$p$ which is endowed with an $O_K$-algebra structure. Denote
by $W(k)$ its ring of Witt vectors. The Frobenius automorphism of $W(k)$ is
denoted by $\sigma$. It is a lifting of the usual Frobenius ${\rm Frob}_k$
of $k$ which sends $x$ to $x^p$. Our main tool will be Dieudonn\'e theory.

\secstart{}\label{\defudm} A {\it Dieudonn\'e module over $k$} is a finitely
generated free $W(k)$-module $M$ together with a $\sigma$-linear
endomorphism $F$ and a $\sigma^{-1}$-linear endomorphism $V$ of $M$ such
that $FV = VF = p\cdot \id_M$. A {\it Dieudonn\'e space over $k$} is a
finite-dimensional $k$-vector space $\Mbar$ together with a ${\rm
Frob}_k$-linear endomorphism $F$ and a ${\rm Frob}^{-1}_k$-linear
endomorphism $V$ such that $FV = VF = 0$. If $(M,F,V)$ is a Dieudonn\'e
module then its reduction modulo $p$ is a Dieudonn\'e space.

A {\it graded Dieudonn\'e module over $k$} is a tuple $((M,F,V), M = M_0 \oplus
M_1)$ where $(M,F,V)$ is a Dieudonn\'e module over $k$ and where $M =
\bigoplus_{i \in \Z/2\Z} M_i$ is a decomposition of $W(k)$-submodules such
that $F$ and $V$ are homogeneous of degree 1. We have the analogous
definition of a {\it graded Dieudonn\'e space over $k$}. If $((\Mbar,F,V),
\Mbar = \Mbar_0 \oplus \Mbar_1)$ is a graded Dieudonn\'e space we call the pair
$$(\dim_k(\Mbar_0/V\Mbar_1), \dim_k(\Mbar_1/V\Mbar_0))$$
its {\it signature} and denote it by ${\rm sign}(\Mbar)$. The {\it
signature} of a graded Dieudonn\'e module is defined as the signature of
its reduction modulo $p$. We have the obvious notion of a {\it morphism of
graded Dieudonn\'e modules} (resp.\ {\it of graded Dieudonn\'e spaces}\/).
Finally by inverting $p$ we get the notion of a {\it graded isocrystal}.

A {\it quasi-unitary Dieudonn\'e module over $k$} is a tuple $\Mline =
((M,F,V), \lrangle, M = M_0 \oplus M_1)$ such that $((M,F,V), M = M_0
\oplus M_1)$ is a graded Dieudonn\'e module over $k$ and $\lrangle\colon M
\times M \ar W(k) \tensor{\Z} \Q$ is a non-degenerate alternating pairing such
that
$$\langle Fx, y \rangle = \langle x, Vy \rangle^{\sigma}$$
and such that $M_0$ and $M_1$ are totally isotropic with respect
$\lrangle$.
We call $\Mline$ a {\it unitary Dieudonn\'e module}\/ if $\lrangle$ is
$W(k)$-valued and perfect.
We have further the obvious notions of a {\it unitary Dieudonn\'e space over
$k$} and a {\it unitary isocrystal over $k$}\/.

\secstart{} By (covariant) Dieudonn\'e theory we have an
equivalence between the category of unitary Dieudonn\'e modules $((M,F,V),
\lrangle, M = M_0 \oplus M_1)$ over $k$ of signature $(r,s)$ and the
category of $p$-divisible $\Iscr$-modules $(X,\iota,\lambda)$ over $k$. Via
this equivalence we have a canonical identification of $\Lie(X)$ with
$M/VM$ such that $L_0 = M_0/VM_1$ and $L_1 = M_1/VM_0$ where $\Lie(X) = L_0
\oplus L_1$ is the decomposition of \defpdivmodule.

\secstart{}\label{\homunitdm} Let $\Mlinebar_1$ and $\Mlinebar_2$ be
two graded Dieudonn\'e spaces over $k$. For every $k$-algebra $R$ we
denote by $\Homline_{\rm gD}(\Mlinebar_1,\Mlinebar_2)(R)$ the set of
$R$-linear maps $\varphi\colon M_1 \tensor{k} R \ar M_2 \tensor{k} R$ which
preserve the grading and such that the diagrams
$$\matrix{M_{1,R} & \arrover{\varphi} & M_{2,R} \cr
\addleft{V} && \addright{V} \cr
M_{1,R}^{(p)} & \arrover{\varphi^{(p)}} & M_{2,R}^{(p)} \cr},\qquad\qquad
\matrix{M_{1,R} & \arrover{\varphi} & M_{2,R} \cr
\auuleft{F} && \auuright{F} \cr
M_{1,R}^{(p)} & \arrover{\varphi^{(p)}} & M_{2,R}^{(p)} \cr}$$
commute. This defines a functor which is representable by a closed
subscheme of $\Homline_k(M_1,M_2)$. In particular it is affine and of
finite type over $k$. If we consider for $\Mlinebar_1 = \Mlinebar_2 =
\Mlinebar$ only those $\varphi$ which are invertible, we get the algebraic
group scheme $\Autline_{\rm gD}(\Mlinebar)$ which is affine as a closed
subgroup scheme of $GL(M)$.

If $\Mlinebar$ is a unitary Dieudonn\'e space over $k$ we denote by
$\Autline_{\rm uD}(\Mlinebar)$ the closed subgroup scheme of automorphisms
which also preserve the symplectic form.


\paragraph{Unitary Dieudonn\'e spaces of signature (n-1,1)}

\secstart{}\label{\possibleisocrystals} From now on $\Mline = ((M,F,V), M =
M_0 \oplus M_1, \lrangle)$ denotes a unitary Dieudonn\'e module over an
algebraically closed field $k$ of characteristic $p$
and we denote by $\Mlinebar$ the unitary Dieudonn\'e space obtained from
$\Mline$ by reduction modulo $p$. We assume that the signature of
$\Mline$ is $(n-1,1)$

Let $N = (M_{\Q}, F \otimes \id)$ be its
isocrystal. For $\lambda \in \Q$ denote by $N_{\lambda}$ ``the'' simple
isocrystal of slope $\lambda$. Then
$$N \cong N(r) \oplus (N_{1/2})^{n-2r}$$
for some integer $0 \leq r \leq n/2$ where
$$N(r) = \cases{0,& if $r = 0$, \cr
N_{{1 \over 2} - {1 \over 2r}} \oplus N_{{1 \over 2} + {1 \over 2r}},& if
$r > 0$ is even, \cr
N_{{1 \over 2} - {1 \over 2r}}^2 \oplus N_{{1 \over 2} + {1 \over 2r}}^2,&
if $r$ is odd. \cr}$$ 
Further the decomposition $N(r) \oplus (N_{1 \over 2})^{n-2r}$ is a
decomposition of unitary isocrystals, i.e.\ it is orthogonal and graded.

It follows that the isocrystal of a unitary Dieudonn\'e module of signature
$(n-1,1)$ is determined by its first Newton slope. If we denote by $p_r$
the Newton-polygon associated to $N(r) \oplus (N_{1/2})^{n-2r}$, we
have $p_0 \leq p_1 \leq \ldots \leq p_{[n/2]}$ where we write $p \leq q$ if
$p$ lies {\it above}\/ $q$.

\secstart{}\label{\basicunitdm} The following two examples will be of vital
importance:
\assertionlist
\assertionitem We define a superspecial unitary Dieudonn\'e module $\Sline$ as
follows: As $W(k)$-module it is a module of rank $2$ generated by elements
$g$ and $h$.
We set $M_0 = W(k)g$ and $M_1 = W(k)h$. The alternating form is given by
$\langle g,h \rangle = 1$, and we have $Fg = -h$ and $Vg = h$. Then $S$ is a
unitary Dieudonn\'e module of signature $(1,0)$. We also write $\Slinebar$
for its reduction mod $p$ and $S$ for the underlying Dieudonn\'e module of
$\Sline$.
\assertionitem For an integer $n \geq 1$ define a unitary Dieudonn\'e module
$\Bline(n)$ as follows: As $W(k)$-module it has a basis
$(e_1,\ldots,e_n,f_1,\ldots,f_n)$ and the grading is given by $M_0 =
W(k)e_1 \oplus \ldots \oplus W(k)e_n$ and $M_1 = W(k)f_1 \oplus \ldots
\oplus W(k)f_n$. The alternating form is defined by
$$\langle e_i,f_j \rangle = (-1)^{i-1}\delta_{ij}.$$
Finally $F$ and $V$ are given by
$$\eqalignno{Ve_i &= f_{i+1},&\hbox{for $i=1,\ldots,n-1$,} \cr
Vf_n &= e_1, \cr
Fe_i &= f_{i-1},&\hbox{for $i = 2,\ldots,n$,} \cr
Ff_1 &= (-1)^ne_n. \cr}$$
This is a unitary Dieudonn\'e module of signature $(n-1,1)$ which we call {\it
braid of length $n$}. We also write $\Blinebar(n)$ for its reduction modulo $p$
and denote the images of $e_i$ and $f_i$ again by $e_i$ resp.\ $f_i$.

\secstart{Lemma}:\label{\braidslope} {\sl The isocrystal associated to
$\Bline(n)$ has the following Newton slopes:
\indention{(2)}
\litem{(1)} If $n$ is odd, every Newton slope is $1/2$.
\litem{(2)} If $n$ is even, the Newton slopes are ${n/2 - 1 \over n}$ and
${n/2 + 1 \over n}$. In particular, it has only slopes $0$ and $1$ if and
only if $n = 2$.}

\proof: Denote by $\lambda(n)$ the first Newton slope of $\Bline(n)$. By
[Zi1] \pz 6.12 we have
$$\lambda(n) = \lim_{n \ar \infty} {1 \over n}\max\set{k \in \Z}{$V^nM
\subset p^kM$}$$
and hence $\lambda(n) = 1/2$ if $n$ is odd and
$\lambda(n) = {n/2-1 \over n}$ if $n$ is even. Therefore the lemma follows from
\possibleisocrystals.

\secstart{Lemma}:\label{\signaturezero} {\sl Assume that the signature of
$\Mline$ is either $(r,0)$ or $(0,r)$ for some $r \geq 1$. Then $\Mline$ is
superspecial, i.e.\ $FM = VM$.}

\proof: Assume that the signature is $(r,0)$ (the other case is
analogous). This means that $VM_0 = M_1$. This implies $F^{-1}M_0^{\bot} =
M_1^{\bot}$ and as we have $M_i^{\bot} = M_{i+1}$ we get $FM_0 = M_1$. It
follows that $FM_1 = VM_1 = pM_0$.

\secstart{Proposition}:\label{\unitdieudsp} {\sl Let $\Mlinebar$ be a
unitary Dieudonn\'e space over $k$ of signature $(n-1,1)$. Then $\Mlinebar$
is isomorphic to
$$\Blinebar(r) \oplus \Slinebar^{n-r}$$
for some integer $r$ with $1 \leq r \leq n$.}

\proof: This follows easily from a theorem of Moonen [Mo]: He has shown (in
much greater generality) that the number of isomorphism classes of unitary
Dieudonn\'e spaces of signature $(n-1,1)$ is given by the number of elements
in the quotient of the Weyl group $\Omega$ of $G$ by the Weyl group
$\Omega_{\mu}$ of the centralizer of $\mu$ in $G_{\C}$. Now $\Omega$ equals
the symmetric group $S_n$ and $\Omega_{\mu}$ can be identified with $S_{n-1}
\times S_1$, embedded in the obvious way in $S_n$. Hence we see that there
are precisely $n$ of such isomorphism classes. As the unitary Dieudonn\'e
spaces $\Blinebar(r) \oplus \Slinebar^{n-r}$ are pairwise nonisomorphic for
different $r$, the claim follows.

\secstart{} In order to study the stratification of the special fibre of
$\Mscr_{\Iscr,K^p}$ by isomorphism types of unitary Dieudonn\'e spaces (the
Oort stratification in the terminology of [We3]) we have to compute the
dimension of the group of automorphisms of a unitary Dieudonn\'e
space. This will be done in the next sections. We refer to [We3] for the
definition of a Dieudonn\'e space over an arbitrary $k$-algebra $R$. We have
the obvious definition of $\Blinebar(n)$ and $\Slinebar$ over $R$.

\secstart{} Let $\Blinebar(n) =
(\Bbar(n), F,V, \lrangle, \Bbar(n) = \Bbar(n)_0 \oplus \Bbar(n)_1)$ be a
braid of length $n$. We use the notations of \basicunitdm(2) and denote the
image of $e_i$ resp.\ $f_i$ in $\Bbar(n)$ again by $e_i$ resp.\
$f_i$. For $r = 0,\ldots,[n/2]$ we set
$$M^r = \bigoplus_{i=0}^{[(n-1)/2]} ke_{n-2i} \oplus \bigoplus_{i =
r}^{[n/2]-1} ke_{n-1-2i}.$$
We can make this definition over arbitrary $k$-algebras. If we write $V^*W
= V^{-1}W \cap \Bbar(n)_{i+1}$ for a subspace $W$ of $\Bbar(n)_i$, we have
$$M^r = (V^*F)^r\Bbar(n)_0.$$
Further, we set
$$t = \cases{n/2,&if $n$ is even,\cr
(n-3)/2,&if $n$ is odd.\cr}$$

\secstart{} Let $\varphi$ be an automorphism of $\Blinebar(n) = (\Bbar(n), F,V,
\lrangle, \Bbar(n) = \Bbar(n)_0 \oplus \Bbar(n)_1)$ over a $k$-algebra
$R$. We claim that $\varphi$ is uniquely determined by $\varphi(e_1)$ and
$\varphi(f_1)$: It is clear that if we know
$\varphi(e_1),\ldots,\varphi(e_i)$ and $\varphi(f_1),\ldots,\varphi(f_i)$
for $i < n$ we have $\varphi(f_{i+1}) = V\varphi(e_i)$. Further
$\varphi(e_{i+1})$ has to satisfy $F(\varphi(e_{i+1})) = \varphi(f_i)$ and
$\langle \varphi(e_{i+1}), \varphi(f_1)\rangle = 0$ and these conditions
determine $\varphi(e_{i+1})$ uniquely as $\varphi(e_{i+1}) \in \Bbar(n)_0$.

As $\varphi(V\Bbar(n)_1) \subset V\Bbar(n)_1$, we have $\varphi(e_1) =
\alpha_{\varphi}e_1$ for a unit $\alpha_{\varphi}$. This defines a
homomorphism $\alpha\colon \Autline(\Blinebar(n)) \arr \G_m$. Set $r =
p^n-1$ if $n$ is even, and $r = p^n+1$ if $n$ is odd. It is easy to see
that the image of $\alpha$ is $\mgd_r$, the scheme of $r$-th roots of unity.

Further, $\varphi(f_1)
\in FM^t$ as $f_1 \in FM^t$. As $\langle e_1,f_i \rangle = \delta_{1i}$,
there exists a unique $m_{\varphi} \in F(M^t)^{(p)}$ such that
$\varphi(f_1) = \alpha_{\varphi}^{-1}(f_1 + m_{\varphi})$ and $\langle
e_1, m_{\varphi}\rangle = 0$, and this means that $m_{\varphi}$ is in the
free submodule generated by $f_2,f_4,\ldots,f_{n-1}$ if $n$ is odd and in
the submodule generated by $f_3,f_5,\ldots,f_{n-1}$ if $n$ is even. Denote
this free submodule by $N$. We consider $N$ as $k$-scheme isomorphic to an
affine space. From the construction above of $\varphi(e_i)$ and
$\varphi(f_i)$ from $\varphi(e_1)$ and $\varphi(f_1)$ it follows that the
morphism of schemes
$$m\colon \Autline(\Blinebar(n)) \asr N, \quad \varphi \asr m_{\varphi}$$
is an isomorphism of schemes (but not of group
schemes). Hence we see:

\claim Proposition: Set $r = p^n-1$ if $n$ is even and $r
= p^n+1$ if $n$ is odd. Then the map
$$\Autline(\Blinebar(n)) \arr
\mgd_r, \qquad \varphi \asr \alpha_{\varphi}$$
provides a surjective homomorphism
of algebraic groups. Its kernel is the identity component
$\Autline(\Blinebar(n))^0$. The morphism
$$m\colon \Aut(\Blinebar(n))^0 \arr N$$
is an isomorphism of varieties.

\secstart{Corollary}: {\sl We have
$$\dim(\Autline(\Blinebar(n))) = \cases{{n \over 2} - 1, &if $n$ is even,\cr
{n - 1 \over 2}, &if $n$ is odd.\cr}$$}

\secstart{Proposition}:\label{\homsspecbraid} {\sl Let $r$ and $m$ be
positive integers.
\assertionlist
\assertionitem The group scheme $\Autline_{\rm uD}(\Slinebar^m)$ is
constant with underlying group equal to the $\F_p$-valued points of the
unitary group in $m$ variables (regarded as an abstract group).
\assertionitem The scheme $\Homline_{\rm gD}(\Slinebar^m,\Blinebar(r))$ is
isomorphic to the $m$-dimensional affine space if $r$ is odd and to
$\alpha_{p^r}^m$ if $r$ is even.}

\proof: Write $\Slinebar^m = \Slinebar_1 \oplus \ldots \oplus \Slinebar_m$
with $\Slinebar_i = \Slinebar$ and for each $1 \geq i \geq m$ choose a
basis $(g_i,h_i)$ as in \basicunitdm(1). Let $R$ be the coordinate
ring of the affine scheme $\Autline_{\rm uD}(\Slinebar^m)$. Let also
$\Phi \in \Autline_{\rm uD}(\Slinebar^m)(R)$ be the universal automorphism. It
preserves the grading, so write $\Phi(g_j)=\sum_ia_{ij}g_i$ and
$\Phi(h_j)=\sum_ib_{ij}h_i$ and denote by $A$ and $B$ the matrices $(a_{ij})$
and $(b_{ij})$. As $\Phi$ commutes with $F$ and $V$, the matrices $A$
and $B$ satisfy
$$A^{(p)} = B, \qquad B^{(p)}= A\leqno{{\rm (A)}}$$
and in particular $A = A^{(p^2)}$, i.e.\ $a_{ij} \in \F_{p^2}$. Moreover,
as $\Phi$ preserves the form, one has
$$\sum_{k=1}^ma_{ki}b_{kj}=\delta_{ij},$$
or equivalently ${}^t\!AB=1$. Therefore $R$ is the quotient of the polynomial
ring $k[a_{ij},b_{ij}]$ by the ideal given by the relations (A) and ${}^t\!AB =
1_m$ which proves (1).

Choose a basis $(e_1,\ldots,e_r, f_1,\ldots,f_r)$ as in \basicunitdm(1) for
$\Blinebar(r)$. To prove (2) we consider the coordinate ring $R$ of
$\Homline_{\rm gD}(\Slinebar^m, \Blinebar(r))$ together with the universal
homomorphism $\Phi \in \Homline_{\rm gD}(\Slinebar^m, \Blinebar(r))(R)$
sending $g_j$ to $\sum_is_{ij}e_i$ and $h_j$ to $\sum_it_{ij}f_i$. The fact
that $\Phi$ commutes with $F$ and $V$ is equivalent to the equations
$$\displaylines{
\rlap{(B)} \hfill t^p_{1,1} = \cdots = t^p_{1,m} = 0, \hfill \cr
\rlap{(C)} \hfill t_{r,1} = \cdots = t_{r,m} = 0, \hfill \cr
\rlap{(D)} \hfill \hbox{$s^p_{i+1,j} = -t_{i,j}$ and $t^p_{i+1,j} =
s_{ij}$} \hfill \cr}$$
for $1 \leq i \leq r-1$ and $1 \leq j \leq m$. The equations (C) and (D)
mean that the only nonzero $s_{ij}$, (resp.\ $t_{ij}$) have $i\equiv r(2)$
(resp.\ $i\equiv r-1(2)$) and, moreover, these are completely determined by
$s_{r,j}$, as $(-1)^{r-i \over 2}s_{r,j}^{p^{r-i}}=s_{i,j}$ if $i\equiv
r(2)$, and $(-1)^{r-i+1 \over 2}s_{r,j}^{p^{r-i}}=t_{i,j}$ if $i\equiv
r-1(2)$.

We have not yet used equation (B), which is non-redundant only if $1\equiv
r-1(2)$, i.e.\ if $r$ is even. In this case it forces $s_{r,j}$ to satisfy
$s_{r,j}^{p^r}=((-1)^{r \over 2}t_{1,j})^p=0$. Therefore 
$$R\cong k[s_{r,1},\ldots,s_{r,m}],$$
if $r$ is odd and
$$R\cong k[s_{r,1},\ldots,s_{r,m}]/(s_{r,1}^{p^r},\ldots,s_{r,m}^{p^r}),$$
if $r$ is even.

\secstart{}\label{\braidplus} We keep the notations. The proof of
\homsspecbraid(2) shows that we have
$$\Homline_{\rm gD}(\Slinebar^m, \Blinebar(r)) = \Homline_{\rm
gD}(\Slinebar^m, \Blinebar^+(r))$$
where we denote by $\Blinebar^+(r)$ the graded Dieudonn\'e subspace which is
generated by
$$e_r,e_{r-2},\ldots,e_{r\ {\rm mod}\ 2}, f_{r-1}, f_{r-3},
\ldots, f_{(r-1)\ {\rm mod}\ 2}.$$
By definition, $\Blinebar^+(r)$ is a totally
isotropic subspace of $\Blinebar(r)$ which is maximal totally isotropic
because of its dimension.

\secstart{Proposition}: {\sl There exists a (non-canonical) isomorphism of
schemes
$$\pi\colon \Autline_{\rm uD}(\Blinebar(r) \oplus \Slinebar^m) \cong
\Homline_{\rm gD}(\Slinebar^m, \Blinebar(r)) \times \Autline_{\rm
uD}(\Slinebar^m) \times \Autline_{\rm uD}(\Blinebar(r)).$$}

\proof: Consider an $R$-valued
point $\varphi$ of $\Homline_{\rm gD}(\Slinebar^m, \Blinebar(r))$ for an
arbitrary $k$-algebra $R$. The restriction of $\varphi$ to $\Slinebar^m
\otimes R$ factors by $\braidplus$ through
$$\varphi'\colon \Slinebar^m \otimes R \arr (\Blinebar^+(r) \oplus
\Slinebar^m) \otimes R.$$
As $\Blinebar^+(r)$ is the radical of $\Blinebar^+(r) \oplus \Slinebar^m$,
one obtains two components
$$\eqalign{\pi_1(\varphi) &\in \Homline_{\rm gD}(\Slinebar^m,
\Blinebar^+(r))(R) = \Homline_{\rm gD}(\Slinebar^m, \Blinebar(r))(R),\cr
\pi_2(\varphi) &\in \Autline_{\rm uD}(\Slinebar^m)(R).\cr}$$
Observe that $\Autline_{\rm uD}(\Slinebar^m)$ is embedded in a canonical
way in $\Autline_{\rm uD}(\Blinebar(r) \oplus \Slinebar^m)$ and that it
acts in a canonical way in $\Homline_{\rm gD}(\Slinebar^m, \Blinebar(r))$
from the right. The maps $\pi_1$ and $\pi_2$ are $\Autline_{\rm
uD}(\Slinebar^m)$-equivariant morphism of schemes. This allows to construct
a section $\sigma$ of the morphism
$$\pi_1 \times \pi_2\colon \Autline_{\rm uD}(\Blinebar(r) \oplus
\Slinebar^m) \arr \Homline_{\rm gD}(\Slinebar^m, \Blinebar(r)) \times
\Autline_{\rm uD}(\Slinebar^m).$$
For this consider an $R$-valued point $(\varphi_1,\varphi_2)$ of the
right hand side. Let $S = (s_{ij})$ (resp. $T = (t_{ij})$) be the ($r
\times m$)-matrix with entries in $R$ which describes $\varphi_1 \circ
\varphi_2^{-1}$ on the bases $g_1,\ldots,g_m$ and $e_1,\ldots,e_r$
(resp. on $h_1,\ldots,h_m$ and $f_1,\ldots,f_r$). Define an
element $\psi$ in $\Autline(\Blinebar(r) \oplus \Slinebar^m)(R)$ as
follows: The $0$-graded component is given by:
$$\eqalign{
e_j &\asr e_j + (-1)^{r-1}\sum_it_{j,i}(g_i+{1 \over 2}\sum_ks_{k,i}e_k), \cr
g_j &\asr g_j+\sum_is_{i,j}e_i, \cr}$$
and the $1$-graded component is given by:
$$\eqalign{
f_j &\asr f_j + (-1)^r\sum_is_{j,i}(h_i+{1 \over 2}\sum_kt_{k,i}f_k), \cr
h_j &\asr h_j + \sum_it_{i,j}f_i. \cr}$$
It is straight forward to see that $\psi$ is compatible with grading, form,
and the operators $F$ and $V$. Further, we have $\pi_1(\psi) = \varphi_1
\circ \varphi_2^{-1}$ and $\pi_2(\psi) = 1$. If we set
$\sigma(\varphi_1,\varphi_2) = \psi \circ \varphi_2$, this defines the
wanted section.

For any $\varphi \in \Autline_{\rm uD}(\Blinebar(r) \oplus \Slinebar^m)(R)$ we
have
$$\varphi\restricted{\Slinebar^m} = \sigma(\pi_1(\varphi),
\pi_2(\varphi))\restricted{\Slinebar^m}.$$
It follows that we have $\pi_3(\varphi) := \varphi^{-1}\sigma(\pi_1(\varphi),
\pi_2(\varphi)) \in \Autline_{\rm uD}(\Blinebar(r))(R)$ where we consider
$\Autline_{\rm uD}(\Blinebar(r))$ as a subgroup scheme of $\Autline_{\rm
uD}(\Blinebar(r) \oplus \Slinebar^m)$ in the canonical way. Then $\pi =
\pi_1 \times \pi_2 \times \pi_3$ is the desired isomorphism.

\secstart{Corollary}:\label{\dimautom} {\sl For $1 \leq r \leq n$ we have
$$\dim{\Autline(\Blinebar(r) \oplus \Slinebar^{n-r})}
= \cases{{r \over 2} - 1,&if $r$ is even,\cr
n - {r + 1 \over 2},&if $r$ is odd.\cr}$$}


\paragraph{Structure of non-supersingular unitary Dieudonn\'e modules}

\secstart{Lemma}\label{\splitawaysspecial}: {\sl Let $\Mline$ be a unitary
Dieudonn\'e module such that its reduction modulo $p$ can be decomposed as
unitary Dieudonn\'e spaces in $\Mlinebar' \oplus \Slinebar^m$ for some $m \geq
1$. Then there exists a decomposition of unitary Dieudonn\'e modules
$\Mline = \Mline' \oplus \Sline^m$ inducing modulo $p$ the given
decomposition.}

\proof: We can assume that $m = 1$. By [Kr] p. 56 we can find a decomposition
of Dieudonn\'e modules $M = M'' \oplus S'$ such that in $S'$ there exists
an $x$ with $Fx = Vx \notin pS'$. The image of $x$ modulo $p$ is a nonzero
element in $(M/pM)_0$. Denote by $x_0$ the $M_0$-component of $x$. Then we
still have $Fx_0 = Vx_0$ and the Dieudonn\'e module $S$ generated by $x_0$
is graded and superspecial. Further the restriction of $\lrangle$ to it is
perfect as it is perfect modulo $p$. This way we get the superspecial
unitary Dieudonn\'e module $\Sline$. It lifts $\Slinebar$. Its orthogonal
complement $M'$ is also a unitary Dieudonn\'e-module $\Mline'$, and it
reduces modulo $p$ to $\Mlinebar'$.

\secstart{}\label{\strucunitdm} Let $\Mline$ be a unitary Dieudonn\'e
module of signature $(n-1,1)$. Denote by $N$ its isocrystal. By
\unitdieudsp\ the reduction modulo $p$ of $\Mline$ is isomorphic to
$\Blinebar(r) \oplus \Slinebar^{n-r}$ for some $r$ with $1 \leq r \leq n$.

\claim Proposition: (1) $N$ is supersingular if and only if $r$ is odd.
\indention{(2)}
\litem{(2)} Assume that $N$ is not supersingular. Consider the
decomposition $N = N' \oplus N''$ with $N' = (N_{1 \over 2})^{n-2m}$ and
$N'' = N(m)$ \possibleisocrystals. Set $M' = M \cap N'$ and $M'' = M \cap
N''$. Then $M = M' \oplus M''$, $M'$ is superspecial, and the
reduction of $M''$ with its induced Dieudonn\'e space structure is
isomorphic to $\Blinebar(2m)$.

\proof: For both assertions we can assume that $r = n$ by
\splitawaysspecial. The condition in (1) is sufficient by \braidslope. Now
let $r > 0$ be even. We show that $N$ has no homogeneous supersingular
subisocrystal. If such a supersingular isocrystal existed we would find a
homogeneous element $x \in N$ of degree 0 such that $Fx = Vx$.
By multiplying with a suitable power of $p$ we can assume that $x
\in M_0 \setminus pM_0$. Modulo $p$ we get a nonzero element $y \in
\Blinebar(r)$ such that $Fy = Vy$. This is a contradiction to
\homsspecbraid(2).

By (1) we know that $r$ is even in the non-supersingular
case and in its proof we showed that $N$ contains no supersingular
subisocrystal. Therefore $M = M''$ (under the assumption $r = n$) and (2)
follows.

\secstart{} Let $\Mline$ be a unitary Dieudonn\'e module
of signature $(n-1,1)$ whose isocrystal is not supersingular. By
\strucunitdm\ we have a decomposition of unitary Dieudonn\'e modules
$\Mline = \Mline' \oplus \Mline''$ where $\Mline'$ is superspecial and
where $\Mlinebar'' \cong \Blinebar(2m)$ (which is equivalent to the fact
that its isocrystal is $N(m)$). The goal of the rest of this chapter is to
give a precise description for those modules. In fact, we will show:

\secstart{Proposition}:\label{\allisbraid} {\sl Let $\Mline$ be a unitary
Dieudonn\'e module whose reduction modulo $p$ is isomorphic to
$\Blinebar(2m)$ for some $m \geq 1$. Then $\Mline \cong \Bline(2m)$.}

\secstart{}\label{\defgenbraid} More generally, for applications to
$p$-isogenies we want to examine the set
$\Mscr$ of quasi-unitary \defudm\ Dieudonn\'e lattices $(M,F,V, M_0 \oplus
M_1)$ of $N(m)$ of signature $(2m-1,1)$ such that
$$M^t := \set{x \in N(m)}{$\langle x,M\rangle \subset W(k)$} = p^{l}M$$
for some $l \in \Z$.

To do this we define the following generalized version of a braid. Let $m
\geq 1$, $l$, and $0 \leq a \leq m-1$ integers. A {\it generalized braid of
lenght $2m$, defect $l$, and jump $a$} is the following graded
Dieudonn\'e module $(M = M_0 \oplus M_1,F,V)$ over $k$ together with a
$W(k)$-bilinear alternating $W(k)_{\Q}$-valued form $\lrangle$:

Let $M_0$ (resp.\ $M_1$) be the free $W(k)$-module with
basis $(e_i)_{i \in \Z/2m\Z}$ (resp.\ with basis $(f_i)_{i \in \Z/2m\Z}$)
such that $F$ is given by
$$\eqalign{Fe_i &= p^{\delta_{i,1}}f_{i-1}, \cr
Ff_i &= p^{1-\delta_{i,2a+1}}e_{i-1}. \cr}$$
Set $V = F^{-1}p$. Finally, define $\lrangle$ by
$$\langle e_i, f_j \rangle = \cases{(-1)^{i-1}p^{l+1}\delta_{ij},&if $i
\leq 2a$,\cr
(-1)^{i-1}p^l\delta_{ij},&if $i > 2a$.\cr}$$
Here and in the sequel we identify $\Z/2m\Z$ with $\{1,\ldots,2m\}$. A
generalized braid is called {\it integral} if $l \geq 0$, i.e.\ $\langle M,
M \rangle \subset W(k)$. It is called {\it quasi-braid}\/ if $a = 0$. We
denote by $\Bscr$ the set of generalized braids in $N(r)$.

\secstart{} From now on let $\Nline = (N,F,V,N_0 \oplus N_1,\lrangle)$ be
a generalized braid of length $2m$, defect $l$ and jump $a$. It follows
from the definition that
$${\rm lg}_{W(k)}(N^t/N) = 4(ml + a).$$
In particular, a generalized braid is a braid iff $a = l = 0$.

If $\Nline$ is a generalized braid, so is $p^n\Nline$ for all $n \in
\Z$. In particular, every $\Mline \in \Mscr$ contains some generalized
braid.

\secstart{Lemma}:\label{\implemma} {\sl Let $\Mline = (M,F,V,M_0 \oplus
M_1,\lrangle)$ be in $\Mscr$ and let $\Nline$ be a generalized braid such
that $N$ is properly contained in $M$. Then we have $V^{-1}N_1 \subset M_0$
or $F^{-1}N_1 \subset M_0$.}

\proof: By replacing the alternating form of $M$ and $N$ by a suitable
$p$-power we can assume that the form is perfect on $M$. As generalized
braids have signature $(2m-1,1)$, one has $\dim V^{-1}N_1/N_0=\dim
F^{-1}N_1/N_0=1$. One can actually say more namely, if $e_1$,...,$e_{2m}$,
$f_1$,...,$f_{2m}$ is a basis of $N$ as in \defgenbraid, then
$$V^{-1}N_1=N_0+W(k)p^{-1}e_{2a},$$
and
$$F^{-1}N_1=N_0+W(k)p^{-1}e_1.$$
If $N$ is properly contained in $M$, then it is a non-perfect, though
integral generalized braid, and we have $<p^{-1}e_{2a},N_1>\ \subset W(k)$
and $<p^{-1}e_1,N_1>\ \subset W(k)$. We conclude that
$$V^{-1}N_1\subset N_1^t, \qquad F^{-1}N_1\subset N_1^t$$
where $N_1^t$ denotes the $W(k)$-lattice in $M_0\otimes\Q$ that is dual to
$N_1$.

Now if $V^{-1}N_1\not\subset M_0$, then the first of the two inclusions
$$M_0\subset V^{-1}N_1+M_0\subset V^{-1}M_1,$$
is proper, so the second one is not because $\dim V^{-1}M_1/M_0=1$.
Consequently $V^{-1}M_1=V^{-1}N_1+M_0\subset N_1^t$ which gives
the inclusion $N_1\subset FM_0$, by considering the dual lattices.

\secstart{} We introduce two maps $\Fscr, \Vscr\colon \Bscr \arr \Bscr$ by
setting
$$\Fscr(N) = N + F^{-1}N_1 + VF^{-1}N_1, \qquad \Vscr(N) = N + V^{-1}N_1 +
FV^{-1}N_1.$$
They satisfy $\Fscr\Vscr = \Vscr\Fscr$, $\Fscr^r\Vscr^r = p^{-1}$, and
$${\rm lg}_{W(k)}(\Fscr(N)/N) = {\rm lg}_{W(k)}(\Vscr(N)/N) = 2.$$

\secstart{Proposition}:\label{\relposbraid} {\sl Let $\Mline = (M,F,V,M_0
\oplus M_1,\lrangle)$ be in $\Mscr$ with $M^t = p^{\lambda}M$ and let
$\Nline$ be a generalized braid of defect $l$ and jump $a$ such that $N \subset
M$. Then we have:
\assertionlist
\assertionitem There exist integers $\alpha, \beta \geq 0$ such that $M
= \Fscr^{\alpha}\Vscr^{\beta}N$.
\assertionitem We have $\alpha + \beta = r(l - \lambda)+a$.
\assertionitem There exists only finitely many $\Mline \in \Mscr$ with $N
\subset M \subset N^t$.
\assertionitem Every $\Mline \in \Mscr$ is a quasi-braid, in particular
\allisbraid\ holds.}

\proof: To prove (1) we make induction on ${\rm
lg}_{W(k)}(M/N)$. If this length is positive, we have either $\Fscr(N)
\subset M$ or $\Vscr(N) \subset M$ by \implemma\ and we can apply the
induction hypothesis to $\Fscr(N) \subset M$ resp.\ $\Vscr(N) \subset M$.

Clearly, (1) implies (4), and (2) follows from the equalities $[M
: N] = [N^t : M^t] = 2(\alpha + \beta)$, $[M^t : M] = 4r\lambda$, and $[N :
N^t] = 4(rl + a)$. Further, (3) follows from (1) and (2): The inclusions $N
\subset M \subset N^t$ implies $\vert[M^t : M]\vert \leq [N^t : N]$, i.e.\
$\vert\lambda\vert \leq l$. With this restriction (and $\alpha, \beta \geq
0$) the equation in (2) has only finitely many solutions.

\secstart{Corollary}:\label{\muorddieudmod} {\sl Let $\Mline$ be a unitary
Dieudonn\`e module of signature $(n-1,1)$ such that the set of slopes of the
associated isocrystal is $\{0, 1/2, 1\}$. Then
$$\Mline \cong \Bline(2) \oplus \Sline^{n-2}.$$}


\paragraph{Newton polygon and Ekedahl-Oort stratification of the special fibre}

\secstart{} From now on we make the assumption $(r,s) = (n-1,1)$ for the
pair $(r,s)$ in \moreno. We denote by $\Mscr$ the special fibre of
$\Mscr_{\Iscr,K^p}$. Let $\F$ be an algebraic closure of $\kappa(O_{E_v})$.

\secstart{Lemma}:\label{\existlemma} {\sl There exists an $\F$-valued point
of $\Mscr$ whose underlying abelian variety is supersingular.}

\medskip

Note that we do not really need this result in the sequel. In fact all
arguments would simplify considerably if the supersingular locus of $\Mscr$
were empty. Hence we remark for the proof only that we can work up to
$\Q$-isogeny. Further it suffices to construct a supersingular abelian
variety $A$ up to isogeny with $B$-action, then we can also find a
polarization on $A$ by [Ko] 9.2. But this can be done using Honda-Tate
theory (see loc.\ cit.\ \pz 10).

\secstart{Theorem}: {\sl Let $s,s'\colon \Spec(k) \arr \Mscr$ be two geometric
points of characteristic $p$ corresponding to tuples
$(A,\iota,\lambda,\hgbar)$ and $(A',\iota',\lambda',\hgbar')$. Assume that
$A$ is not supersingular. Then the following statements are equivalent:
\indention{(2)}
\litem{(1)} The $p$-divisible $\Iscr$-modules
$(A,\iota,\lambda)[p^{\infty}]$ and $(A',\lambda',\iota')[p^{\infty}]$ are
isomorphic.
\litem{(2)} The $BT_1$'s with $O_B$-action and pairing
$(A,\iota,\lambda)[p]$ and $(A',\lambda',\iota')[p]$ are isomorphic.
\litem{(3)} The isocrystals with $B$-action and pairing
associatd to $(A,\iota,\lambda)$ and $(A',\lambda',\iota')$ are isomorphic.}

\proof: In each of the three cases (1) -- (3) the non-supersingularity of
$A$ implies that $A'$ is not supersingular (this is obvious for (1) and (3)
and it follows from \strucunitdm(1) if (2) holds). By \allisbraid\ and
\splitawaysspecial, (1) and (2) are equivalent. Finally \strucunitdm(2)
implies the equivalence of (2) and (3).

\secstart{}\label{\describestrat} Let $\Spec(k) \arr \Mscr$ be a geometric
point of characteristic $p$ corresponding to a tuple
$(A,\iota,\lambda,\hgbar)$. The Dieudonn\'e space of the $p$-torsion gives
a unitary Dieudonn\'e space of signature $(n-1,1)$ which is isomorphic to
$\Blinebar(\rho(\sbar)) \oplus \Slinebar^{n-\rho(\sbar)}$ for some integer
$\rho(\sbar) \in \{1,\ldots,n\}$ by \unitdieudsp. This defines the
Ekedahl-Oort stratification
$$\Mscr = \bigcup_{\rho = 1}^n \Mscr_{\rho}.$$
Now [We3] shows that all $\Mscr_{\rho}$ are equi-dimensional and that if
$\Mscr_{\rho} \not= \emptyset$ we have
$$\codim(\Mscr_{\rho}, \Mscr) = \dim{\Autline(\Blinebar(\rho) \oplus
\Slinebar^{n-\rho})} = \cases{{\rho \over 2} - 1,&if $\rho$ is even,\cr
n - {\rho + 1 \over 2},&if $\rho$ is odd\cr}\lformno\sublabel{\dimOortstrata}$$
by \dimautom. By \existlemma\ and \strucunitdm\ we know that the union of
the $\Mscr_{\rho}$ where $\rho$ runs through the odd numbers between 1 and
$n$ is nonempty as this is precisely the supersingular locus.

In this case the Ekedahl-Oort stratification is finer than the Newton polygon
stratification \strucunitdm\ which is given by
$$\Mscr = \Mscr_2 \cup \Mscr_4 \cup \ldots \cup \Mscr_{2[n/2]} \cup
\bigcup_{\rm \rho\ odd} \Mscr_{\rho}.$$
Further de Jong and Oort have shown in [dJO] that going from one Newton
polygon stratum $\Mscr_{\rho}$ to the next smaller one, the codimension can
only jump by 1. It follows therefore from \dimOortstrata\ that all Newton
polygon strata are non-empty and equidimensional. Further,
$\Mscr_{\rho^{\rm ss}}$ is dense in the supersingular locus
$\bigcup_{\hbox{$\rho$ odd}} \Mscr_{\rho}$ of $\Mscr$ where $\rho^{\rm ss}$ is
the maximal odd integer $\leq n$.

\secstart{Proposition}:\label{\dimsupersing} {\sl The dimension of the
supersingular locus is equal to $[(n-1)/2]$.}

\proof: This follows from \describestrat.

\secstart{} The stratum $\Mscr_2$ of $\Mscr$ is the $\mu$-ordinary stratum
in the sense of [We1]. We also denote it by $\Mscr^{\mu}$. If $\Spec(k)
\arr \Mscr^{\mu}$ is a geometric point corresponding to
$(A,\iota,\lambda,\hgbar)$ with associated unitary Dieudonn\'e module $(M,
F,V, M = M_0 \oplus M_1, \lrangle)$, the associated slope filtration of
$(M,F,V)$ is actually a graduation
$$M = M(0) \oplus M(1/2) \oplus M(1)$$
and the $M(\lambda)$ are homogeneous Dieudonn\'e submodules for $\lambda
= 0, 1/2, 1$; $M(0)$ and $M(1)$ are totally isotropic and orthogonal to
$M(1/2)$ and they are free $O_K \otimes W(k)$-modules of rank 1; $M(1/2)$
is orthogonal to $M(0) \oplus M(1)$ and it is isomorphic to $\Sline^{n-2}$.


\paragraph{A generic model for the moduli space of p-isogenies}

\secstart{} The goal in this chapter is to show that there are no components of
$\pisog \otimes \kappa(E_v)$ \definepisog\ of dimension $< n-1$. For this
we use Grothendieck-Messing theory to translate this problem into linear
algebra.

\secstart{Definition}: Let $S$ be any $\Z_p$-scheme and let $c \geq 0$ be
an integer. A {\it $p$-isogeny
  pair over $S$ of type $\Iscr$ of multiplicator $c$} or shorter an {\it
  $\Iscr$-pair over $S$} is a tuple $((M,\lambda),(M',\lambda'),\alpha)$ where
\bulletlist
\bulletitem $M$ and $M'$ are $O_K \tensor{\Zp} \Oscr_S$-modules which
  are locally on $S$ free $\Oscr_S$-modules of rank $2n$.
\bulletitem $\lambda\colon M \arriso M\vdual$ and $\lambda'\colon M' \arriso
  {M'}\vdual$ are isomorphisms of $O_K \tensor{\Zp} \Oscr_S$-modules such
  that $\lambda\vdual = -\lambda$ and ${\lambda'}\vdual = -\lambda'$. Here
  we denote by $(\ )\vdual$ the $\Oscr_S$-linear dual which is
  endowed with an $O_K$-module structure via $b\cdot m\vdual(m) =
  m\vdual(b^*m)$ for $b \in O_K$, $m \in M$ and $m\vdual \in M\vdual$.
\bulletitem $\alpha\colon M \arr M'$ is an $O_K \tensor{\Zp}
  \Oscr_S$-linear map such that
$$\alpha\vdual \circ \lambda'\circ \alpha = p^c\lambda \quad \hbox{and}
\quad \alpha \circ \lambda^{-1} \circ \alpha\vdual = p^c{\lambda'}^{-1}.$$

A morphism of $\Iscr$-pairs is a commutative diagram
$$\matrix{M_1 & \arrover{\alpha_1} & M'_1 \cr
\addleft{\psi} && \addright{\psi'} \cr
M_2 & \arrover{\alpha_2} & M'_2 \cr}$$
such that $\psi$ and $\psi'$ are $O_B \tensor{\Z_{(p)}} \Oscr_S$-linear
symplectic isomorphisms.

For varying $S$ the $\Iscr$-pairs over $S$ with multiplicator $c$ define an
Artin stack ${\tt dR}_{\Iscr}^c$ over $\Zp$, i.e.\ an algebraic stack in
the sense of [LM].

\secstart{} Let $S$ be a scheme. Denote by $\Ptt^c(S)$ the category of
tuples $(N,N',u,v)$ where $N$ and $N'$ are locally free $\Oscr_S$-modules
of rank $n$ and where
$$u\colon N \arr N', \qquad v\colon N' \arr N$$
are $\Oscr_S$-linear maps such that $v\circ u = p^c\id_N$ and $u
\circ v = p^c\id_{N'}$. A morphism
$$(N_1,N'_1,u_1,v_1) \arr (N_2,N'_2,u_2,v_2)$$
in this category is a pair of isomorphisms $\psi\colon N_1 \arr N_2$ and
$\psi'\colon N'_1 \arr N'_2$ such that the diagrams
$$\matrix{N_1 & \arrover{u_1} & N'_1 \cr
\addleft{\psi} && \addright{\psi'} \cr
N_2 & \arrover{u_2} & N_2 \cr}, \qquad
\matrix{N_1 & \allover{v_1} & N'_1 \cr
\addleft{\psi} && \addright{\psi'} \cr
N_2 & \allover{v_2} & N_2 \cr}$$
commute.

We have the obvious notion of a pullback with respect to a morphism $S'
\arr S$ and this makes $\Ptt^c$ into an algebraic stack over $\Z$.

\secstart{}\label{\reformlindat} Let $S$ be an $O_K$-scheme and let
$((M,\lambda),(M',\lambda'),\alpha)$ be an $\Iscr$-pair over $S$ of
multiplicator $c$. The
decomposition $O_K \tensor{\Zp} \Oscr_S = \Oscr_S \times \Oscr_S$ induced
by \unitdecomp\ gives a decomposition $M = M_0 \oplus M_1$ such that $M_i$
is totally isotropic with repsect to the bilinear form given by
$\lambda$. Hence $\lambda$ defines an isomorphisms $l_0\colon M_0 \arriso
M\vdual_1$ and $l_1\colon M_1 \arriso M\vdual_0$. The same holds for
$(M',\lambda')$. 
The map $\alpha$ is homogeneous with respect to the decomposition, i.e.\
$\alpha = \alpha_0 \oplus \alpha_1$ for $\Oscr_S$-linear maps $\alpha_i\colon
M_i \arr M'_i$. The functor
$$\Fscr\colon ((M,\lambda),(M',\lambda'),\alpha) \asr
(M_0,M'_0,\alpha_0,l_0^{-1} \circ \alpha\vdual_1 \circ l_0')$$
defines an equivalence of the category of $\Iscr$-pairs of
multiplicator $c$ over $S$ with $\Ptt^c_S$.

This equivalence is compatible with pullback via morphisms $S' \arr
S$. Hence we get an isomorphism of the algebraic stacks
$$\Fscr\colon {\tt dR}_{\Iscr}^c \tensor{\Zp} O_K \arriso \Ptt^c \otimes_{\Z}
O_K.$$

\secstart{} We are going to study the special fibre of $\pisog^c$. More
precisely we will study
$$\pisog^c_0 = \pisog^c \otimes \kappa(O_K).$$
For $c = 0$ each one of the morphisms $s,t\colon \pisog^c \arr
\Mscr_{\Iscr,K^p}$ is an isomorphism. Hence we will assume form now on that
$c > 0$.

We will define a stratification of $\pisog^c_0$ by locally closed
subschemes using the algebraic stack
$$\Ptt := \Ptt^c \tensor{\Z} \kappa(O_K).$$
We remark that $\Ptt$ is independent of $c > 0$.

\secstart{}\label{\stratmorph} Let $S$ be a $\kappa(O_K)$-scheme and let
$f\colon (A,\iota,\lambda,\hgbar) \arr (A',\iota',\lambda',\hgbar')$ be an
$S$-valued point of $\pisog^c$. By taking deRham cohomology we get via
functoriality a homomorphism
$$\alpha = H_{DR}(f)\colon M = H_{DR}(A) \arr M' = H_{DR}(A').$$
Further, $\lambda$ and $\iota$ endow $H_{DR}(A)$ with a homomorphisms
$\lambda\colon M \arr M\vdual$ such that $\lambda\vdual = -\lambda$ and an
$O_K$-action; ditto for  $M'$. By the definition we have
$$\alpha\vdual \circ \lambda'\circ \alpha = p^c\lambda = 0.$$
By considering the dual $p$-isogeny $f\vdual$ and using the fact that
$H_{DR}(f\vdual) = -\alpha\vdual$ we also have the relation
$$\alpha \circ \lambda^{-1} \circ \alpha\vdual = p^c{\lambda'}^{-1} = 0.$$
Altogether we get a morphism $\pisog^c_0 \arr {\tt dR}_{\Iscr}^c \otimes
\kappa(O_K)$. By composition with the functor $\Fscr$ \reformlindat\ we
obtain a morphism
$$\Psi\colon \pisog^c_0 \arr \Ptt.$$

\secstart{Proposition}:\label{\describepairs} {\sl Let $S$ be a
$\kappa(O_K)$-scheme and let $\Ntilde_1 = (N_1,N_1',u_1,v_1)$ and
$\Ntilde_2 = (N_2, N'_2, u_2,v_2)$ be two objects of $\Ptt(S)$. Assume
that $\Im(u_i)$ is locally a direct summand of $N'$ of rank $m$ and that
$\Im(v_i)$ is locally a direct summand of $N$ of rank $l$ (where $m$ and
$l$ are nonnegative integers independant of $i = 0,1$). Let
$\Isomline_S(\Ntilde_1,\Ntilde_2)$ be the $S$-scheme of isomorphisms
between $\Ntilde_1$ and $\Ntilde_2$.
\assertionlist
\assertionitem $\Isomline_S(\Ntilde_1,\Ntilde_2)$ admits Zariski locally on
$S$ sections.
\assertionitem $\Isomline_S(\Ntilde_1,\Ntilde_2)$ is smooth of
relative dimension
$$d(m,l) := n^2 + (n - m - l)^2$$
over $S$.}

\proof: To prove (1) it suffices to show that there exists locally on $S$ a
normal form of $\Ntilde_1$ and $\Ntilde_2$ which depends only on $m$,
$l$. To shorten notations we write $(N,N',u,v)$ for $\Ntilde_i$. Hence (1)
follows from the following obvious lemma

\claim Lemma: Let $S = \Spec(R)$ be a local scheme. There exist
isomorphisms $\varphi\colon N \arriso R^n$ and $\varphi'\colon N' \arriso
R^n$ such that $u$ and $v$ are given by the matrices 
$$U = \pmatrix{I_m & 0_{m,n-m} \cr 0_{n-m,m} & 0_{n-m,n-m} \cr}, \qquad
V = \pmatrix{0_{m,m} & 0_{m,l} & 0_{m,n-m-l} \cr
0_{l,m} & I_l & 0_{l,n-m-l} \cr
0 & 0 & 0 \cr}$$
where $0_{i,j}$ denotes the zero matrix of size $i \times j$.

\medskip

\proof of the Proposition (continued): To prove (2) we can assume that
$\Ntilde_1 = \Ntilde_2 =: \Ntilde$ is in the normal form of the lemma. It
is clear that $\Autline_S(\Ntilde)$ is smooth and a straight forward
calculation gives its relative dimension.

\secstart{Lemma}:\label{\Fittinglemma} {\sl Let $S$ be a scheme and let
$\Escr$ be a quasicoherent $\Oscr_S$-module of finite type. For every
integer $n \geq 0$ denote by ${\rm rk}_{= n}(\Escr)$ the functor on
$({\tt Sch}/S)$
$$(f\colon T \arr S) \asr \cases{\{f\},&if $f^*\Escr$ is locally free of
rank $n$,\cr
\emptyset,&otherwise.}$$
The functor ${\rm rk}_{= n}(\Escr)$ is representable by a
locally closed subscheme $Z_n$ of $S$ and the immersion $Z_n \arr S$ is of
finite presentation if $\Escr$ is of finite presentation.}

\proof: We can assume that $S = \Spec(A)$ is affine, $M$ the finitely
generated $A$-module corresponding to $\Escr$. We define $Y_n$ as the
closed subscheme given by the $(n-1)$-th Fitting ideal ${\rm
Fitt}_{n-1}(M)$ of $M$ (see e.g.\ [Eis] \pz 20.2) where we make the convention
${\rm Fitt}_{-1}(M) = 0$. If $M$ is of finite presentation
the construction of ${\rm Fitt}_n(M)$ shows that the Fitting ideals are
finitely generated. Set
$$Z_n = Y_n \setminus Y_{n-1}.$$
This represents the functor ${\rm rk}_{= n}(\Fscr)$ by loc.\ cit.\ 20.5--8.

\secstart{} Let $\Xi = \Xi_n$ be the set of pairs $(m,l)$ of nonnegative
integers such that $m + l \leq n$. For $\xi_0
= (m_0,l_0)$ let $U_{\xi_0}$ be the set of pairs $(m,l) \in \Xi$
such that $m \geq m_0$ and $l \geq l_0$. We endow
$\Xi$ with the topology generated by the $U_{\xi_0}$ for $\xi_0 \in
\Xi$. For $\xi = (m,l) \in \Xi$ we denote by $\Ptt^{\xi}$ the substack of
$\Ptt$ of those $(N,N',u,v)$ where $\Im(u)$ is a
direct summand of rank $m$ and $\Im(v)$ is a direct summand of rank
$l$. Applying \Fittinglemma(2) to the cokernels of $u$ and $v$ we see that
the inclusion $\Ptt^{\xi} \air \Ptt$ is representable by an immersion. The
underlying topological space of $\Ptt^{\xi}$ consists of one point as all
objects in $\Ptt^{\xi}$ are locally isomorphic \describepairs. By
\Fittinglemma(1) we further have (cf.\ [We3] 4.4): 

\claim Corollary: The underlying topological spaces of $\Ptt$ is equal to
$\Xi$.

\secstart{} Let $S$ be a $\kappa(O_K)$-scheme. We denote by $\Mtt(S)$ the
category of tuples $(N,N',u,v,H,H')$ where $(N,N',u,v)$ is
an object in $\Ptt(S)$ and where $H \subset N$ and $H' \subset N'$ are
locally direct summands of rank $1$. A morphism
$$(N_1,N'_1,u_1,v_1,H_1,H_1') \arr (N_2,N_2',u_2,v_2,H_2,H_2')$$
in $\Mtt(S)$ is a morphism $(\psi,\psi')$ in $\Ptt(S)$ such that $\psi(H_1)
\subset H_2$ and $\psi'(H'_1) \subset H'_2$.

For varying $S$ we get the algebraic stack $\Mtt$ over $\kappa(O_K)$.

\secstart{} Associating to a $p$-isogeny the induced morphism on the deRham
cohomology together with the Hodge filtration gives via the isomorphism
$\Fscr$ as in \stratmorph\ a morphism
$$\Phi\colon \pisog_0^c \arr \Mtt.$$
The forgetful functor $(N,N',u,v,H,H') \asr (N,N',u,v)$ defines a morphism
$$\Pi\colon \Mtt \arr \Ptt$$
and we have
$$\Psi = \Pi \circ \Phi.$$
For $\xi \in \Xi$ we set
$$\eqalign{\Mtt^{\xi} &:= \Pi^{-1}(\Ptt^{\xi}),\cr
(\pisog^c_0)^{\xi} &:= \Psi^{-1}(\Ptt^{\xi}).\cr}$$
These is a locally closed substack of $\Mtt$ resp.\ a locally closed
subscheme of $\pisog^c_0$.

\secstart{}\label{\Mstructure} For $\xi = (m,l) \in \Xi$ let $\Ntilde^{\xi}_0 =
(N_0,N'_0,u_0,v_0)$ be the object in $\Ptt^{\xi}$ which has the normal form
as in \describepairs. Denote by $\Mtttilde^{\xi}$ the functor which associates
to an $O_K$-scheme $S$ the set of isomorphism classes of tuples
$(N,N',u,v,H,H',\alpha)$ where $(N,N',u,v,H,H')$ is an object in
$\Mtt^{\xi}$ and where $\alpha$ is an isomorphism $(N,N',u,v) \arriso
\Ntilde^{\xi}_0 \tensor{\kappa(O_K)} \Oscr_S$ is isomorphic to the scheme
whose $S$-valued points are given by pairs $(H,H')$ where $H \subset N_0
\otimes \Oscr_S$ and $H' \subset N'_0 \otimes \Oscr_S$ are locally direct
summands of rank 1 such that $(u_0 \otimes \id)(H) \subset H'$ and $(v_0
\otimes \id)(H') \subset H$. Hence it is representable by a
projective scheme over $\kappa(O_K)$. By \describepairs\ the canonical morphism
$$\Mtttilde^{\xi} \arr \Mtt^{\xi}$$
is a torseur under the smooth $\kappa(O_K)$-group scheme
$\Autline(\Ntilde^{\xi}_0)$.

We consider $\Mtttilde^{\xi}$ as a closed subscheme of $\P(N_0)
\times \P(N_0')$. It has three distinguished closed subschemes $Z$, $Z'$ and
$Z''$: $Z$ is the closure of the locus where $(u_0 \otimes \id)(H)$ is
locally a direct summand of rank 1, and it is
isomorphic to $\P(N_0)$ via the first projection $\Mtttilde^{\xi} \arr
\P(N_0)$; $Z'$ is the closure of the locus where $(v_0 \otimes \id)(H')$ is
locally a direct summand of rank 1, and it is
isomorphic to $\P(N'_0)$ via the second projection $Z \arr \P(N_0')$; $Z''$
is the locus where $(u_0 \otimes \id)(H) = (0)$ and $(v_0 \otimes \id)(H') = (0)$ and it is isomorphic to
$\P(\Ker(u_0)) \times \P(\Ker(v_0))$ via the embedding $Z \air \P(N_0) \times
\P(N_0')$. Further we have $Z \cap Z'' = \P(\Ker(u_0)) \times \P(\Im(u_0))
\subset Z''$ and $Z' \cap Z'' = \P(\Im(v_0)) \times \P(\Ker(v_0))
\subset Z''$. Hence we get:

\claim Lemma: If $m = n$ or $l = n$, we have $\Mtttilde^{\xi} \cong
\P^n$. If $m + l = n$ and $m,l \geq 1$, $\Mtttilde^{\xi}$ has the two
components $Z$ and $Z'$ and their intersection has codimension 1 in $Z$ and
$Z'$. If $m + l < n$ we have three components, namely $Z$, $Z'$ and
$Z''$. In any case we have $\dim(\Mtilde^{\xi}) \geq n-1$.

\secstart{Proposition}:\label{\locmodpisog} {\sl Fix a $\xi = (m,l)$ in
$\Xi$. Let $U$ be an open subscheme of $\pisog^c_0$ which is contained in
$(\pisog^c_0)^{\xi}$. The morphism
$$\Phi_U\colon U \arr \Mtt^{\xi}$$
induced by $\Phi$ is smooth of relative dimension $d(m,l)$ \describepairs.}

\proof: We first prove smoothness using Grothendieck-Messing theory. As
$\Phi_U$ is a morphism of finite type between noetherian Artin stacks it
suffices to show that for all local Artinian rings $R$ and for all
quotients $R_0$ of $R$ defined by an ideal $I$ of square zero and for all
commutative diagrams
$$\matrix{\Spec(R_0) & \arrover{\alpha} & U\cr
\add && \addright{\Phi_U}\cr
\Spec(R) & \arrover{\beta} & \Mtt^{\xi}\cr}$$
there exists a morphism $\Spec(R) \arr U$ which commutes with the
diagram. As $U$ is open in $\pisog^c_0$ it suffices to construct a morphism
$\Spec(A) \arr \pisog^c_0$ which makes the corresponding diagram
commutative. As $I^2 = (0)$ we can equip $I$ with a PD-structure. The
morphism $\alpha\colon \Spec(A_0) \arr U$ corresponds to a $p$-isogeny $f\colon
(A_1,\iota_1,\lambda_1,\hgbar_1) \arr (A_2,\iota_2,\lambda_2,\hgbar_2)$ over
$R_0$, taking the corresponding morphisms of crystals and evaluating it in
$R$ we get an $R$-valued point of $\Ptt^{\xi}$. By \describepairs\ we find an
isomorphism of this $R$-valued point of $\Ptt^{\xi}$ with the $R$-valued point
given by the composition of $\beta$ with $\Mtt^{\xi} \arr \Ptt^{\xi}$
induced by $\Pi$. Via this isomorphism the morphism $\beta$ defines a
lifting of the Hodge filtration in the crystal of $f$ hence by
Grothendieck-Messing theory the desired morphism $\Spec(R) \arr \pisog^c_0$.

To prove the assertion about the relative dimension let $\Utilde$ be the
fibre product of $U$ and
$\Mtttilde^{\xi}$ over $\Mtt^{\xi}$. Then $\Utilde \arr U$ is an
$\Autline(\Ntilde^{\xi}_0)$-torseur and hence smooth of relative dimension
$d(\xi)$ \describepairs. In order to show that the smooth morphism $\Utilde
\arr \Mtttilde^{\xi}$ has relative dimension $d(\xi)$, it suffices to prove
that if $\utilde$ is some $k$-valued point of $\Utilde$ with image $u$ resp.\
$\mtilde$ in $U$ resp.\ $\Mtttilde^{\xi}$ we have $\dim_kT_u(U) =
\dim_kT_{\mtilde}\Mtttilde^{\xi}$ which again follows as above from
\describepairs\ and Grothendieck-Messing theory applied to $R_0 = k$ and
$R = k[\eps]/\eps^2$.

\secstart{Corollary}:\label{\locstructure} {\sl We keep the notations of
the proof of \locmodpisog. For all closed points $u \in U$ there exists an
\'etale neighborhood $U'$ of $u$ and a section $U' \arr \Utilde$ such that the
composition $U' \arr \Utilde \arr \Mtttilde^{\xi}$ is \'etale. In
particular, for all closed points $u \in U$ there exists an $\mtilde \in
\Mtttilde^{\xi}$ such that $\dim_u(U) = \dim_{\mtilde}\Mtttilde^{\xi}$.}

\proof: As we need only the second assertion, we omit the proof of the more
precise statement. To prove the second assertion we use the notations of
the proof of \locmodpisog. We have seen that there are surjective smooth
morphism $\Utilde \arr U$ and $\Utilde \arr \Mtttilde^{\xi}$ of the same
relative dimension $d(\xi)$ which implies the claim.

\secstart{Proposition}:\label{\dimprop} {\sl Every component of $\pisog^c
\otimes \kappa(O_{E,\nu})$ has dimension $\geq n-1$.}

\proof: We can assume that $c > 0$. It suffices to show that every
component of $\pisog^c_0$ has dimension $\geq n-1$. Let
$X$ be a component of $\pisog^c_0$ and let $\xi$ be the
element of $\Xi$ such that the generic point $\eta$ of $X$ lies in
$(\pisog^c_0)^{\xi}$. Hence $(\pisog^c_0)^{\xi} \cap X$ is open in $X$. If
$U$ is an open subset of $(\pisog^c_0)^{\xi} \cap X$ which does not meet
any other component of $\pisog^c_0$, $U$ is open in $\pisog^c_0$ and dense
in $X$. We can apply \locstructure\ to $U$ and the claim follows from the
fact that every component of $\Mtttilde^{\xi}$ has dimension $\geq n-1$
\Mstructure.


\paragraph{Source and target morphism}

\secstart{}\label{\notatfibre} To shorten notations we set $Y = Y^c = \pisog^c
\otimes \kappa(O_{E_v})$. Let $Y^{\mu} \subset Y$ be the $\mu$-ordinary
locus and denote by $Y^{\rm nss}$ (resp.\ $\Mscr^{\rm nss}$) the
non-supersingular locus of $Y$ (resp.\ of the special fibre $\Mscr$ of
$\Mscr_{\Iscr,K^p}$). Let $k$ be an algebraic closure of $\kappa(O_K)$ and
let $\xbar = (A,\iota,\lambda,\hgbar)$ be a $k$-valued point of $\Mscr^{\rm
nss}$. Denote by $Y^{\rm nss}_{\xbar}$ the fibre of the source morphism
$Y^{\rm nss} \arr \Mscr^{\rm nss}$ in $\xbar$.

\secstart{Proposition}:\label{\finitefibres} {\sl The number of $k$-valued
points in $Y^{\rm nss}_{\xbar}$ is finite and depends only on the
isogeny class of $\xbar$.}

\proof: Let $\Mline = (M,F,V, M = M_0 \oplus M_1, \lrangle)$ be a unitary
Dieudonn\'e module associated to $(A,\iota,\lambda,\hgbar)$. To give a
$k$-valued point of $Y_{\xbar}$ is the same as to give an equivalence
class of homomorphisms
$$\alpha\colon \Mline \arr \Mline'\leqno{(*)}$$
of graded Dieudonn\'e modules such that
$\langle \alpha(x), \alpha(y)\rangle' = p^c\langle x,y \rangle$ where
$\Mline'$ is another unitary Dieudonn\'e module. Here we call
$(\alpha_1,\Mline'_1)$ and $(\alpha_2,\Mline'_2)$ equivalent if there
exists an isomorphism $\psi\colon \Mline'_1 \arriso \Mline'_2$ such that
$\psi \circ \alpha_1 = \alpha_2$. We have decompositions $\Mline = B \oplus
S$ and $\Mline' = B' \oplus S'$ of unitary Dieudonn\'e modules where $B
\cong B' \cong \Bline(r)$ and $S \cong S' \cong \Sline^{n-r}$ for
some even integer $2 \leq r \leq n$ \splitawaysspecial\ and
\allisbraid. As the isogeny type depends only on $r$ \possibleisocrystals,
the second claim follows.

The homomorphism $\alpha$ respects
this decomposition. Hence to show that there are only finitely many
$\alpha$ up to equivalence as in (*) we can assume that either $\Mline
\cong \Mline' \cong \Bline(r)$ or that $\Mline \cong \Mline' \cong
\Sline^{n-r}$. The claim for the first case follows from
\relposbraid. Hence we can assume that we are in the second case. It
suffices to show that for any integer $m \geq 1$ there are only finitely
many equivalence classes of homomorphisms of graded Dieudonn\'e modules
$\beta\colon \Sline^m \arr \Sline^m$ such that the $W(k)$-length of the
cokernel of $\beta$ is equal to $cm/2$. But any homomorphism
$\beta$ is determined by $f := \beta\restricted{S_0^m}\colon S_0^m \arr
S_0^m$. As $\beta$ commutes with $F$ and $V$ we have $f^{\sigma^2} = f$,
i.e. $f \in M_m(W(\F_{p^2}))$. The compatibility with the alternating forms
translates into the condition ${}^tf^{\sigma}f = p^c$. Two such elements
$f_1$ and $f_2$ are equivalent if and only if there exists an $h \in
J(\Z_p)$ where $J$ is the $\Zp$-group scheme
$$J = \set{g \in {\rm Res}_{W(\F_{p^2})/\Z_p}GL_m}{${}^tg^{\sigma}g = 1$}$$
such that $hf_1 = f_2$. Hence we have to show that the $J$-orbits on
$$X_c = \set{f \in M_m(W(\F_{p^2}))}{${}^t\!f^{\sigma}f = p^c$}$$
are finite. Every $J(\Zp)$-orbit on $X_c$ is open as $J(\Zp)$ is open in
$J(\Qp)$ and for any two points $x,x' \in X_c$ there exists a $g \in
J(\Qp)$ such that $gx = x'$. As $X_c$ is compact, $X_c/J(\Z_p)$
is compact and discrete and therefore finite.

\secstart{Corollary}:\label{\sourcefinite} {\sl The restriction of source and
target morphism to $Y^{\rm nss}$ is a finite morphism $Y^{\rm nss} \arr
\Mscr^{\rm nss}$.}

\proof: The source morphism is quasi-finite by \finitefibres, hence it is
finite because it is also proper. By duality the same holds for the target
morphism.

\secstart{}\label{\expodef} Let $f\colon (A,\iota,\lambda,\hgbar) \arr
(A',\iota',\lambda',\hgbar')$ be a $k$-valued point of $\pisog^{\mu}$ where
$k$ is an algebraically closed field. Let $\Cscr$ be the category of pairs
$(R,\alpha)$ where $R$ is a local Artinian $O_{E_v}$-algebra and where
$\alpha$ is an isomorphism $\kappa(R) \arriso k$. We denote by $\Def(f)$
the functor of deformations of $f$ to $\Cscr$. To give a deformation of $f$
is the same as to give a deformation of the morphism
$$f\colon \Xline = (X,\iota,\lambda) \arr \Xline' = (X',\iota',\lambda')$$
of $p$-divisible $\Iscr$-modules associated to $f$ \assocdm\ by the theorem
of Serre-Tate. 

There exists a unique decomposition $X = X_{\inf} \times X_{\et}$ where
$X_{\inf}$ and $X_{\et}$ are the infinitesimal resp.\ the \'etale part of
$X$ with induced $O_K$-action $\iota_{\inf}$ and $\iota_{\et}$. Further
$\lambda$ induces isomorphisms $\lambda_{\inf}\colon X_{\inf} \arriso
(X\vdual)_{\inf}$ and $\lambda_{\et}\colon X_{\et} \arriso
(X\vdual)_{\et}$. We set $\Xline_{\inf} =
(X_{\inf},\iota_{\inf},\lambda_{\inf})$, $\Xline_{\et} =
(X_{\et},\iota_{\et},\lambda_{\et})$ and $\Xline_{\bi} =
((\Xline_{\inf}\vdual)_{\inf})\vdual$.
We make the same definitions of $\Xline'$ and get induced morphisms
$f_{\inf}$, $f_{\et}$ and $f_{\bi}$.

We have canonical morphisms $\Def(f) \arr \Def(f_{\inf}) \arr
\Def(f_{\bi})$. Further, sending a deformation to its source (resp.\
target) gives morphisms $\Def(f) \arr \Def(\Xline)$, $\Def(f_{\inf}) \arr
\Def(\Xline_{\inf})$ and $\Def(f_{\bi}) \arr \Def(\Xline_{\bi})$
(resp. $\Def(f) \arr \Def(\Xline')$, $\Def(f_{\inf}) \arr
\Def(\Xline'_{\inf})$ and $\Def(f_{\bi}) \arr \Def(\Xline'_{\bi})$).

\secstart{Proposition}:\label{\describedef} {\sl With the notations of
\expodef\ there are the following relations and descriptions of these functors:
\assertionlist
\assertionitem The induced morphisms
$$\eqalign{\Def(f) &\arr \Def(\Xline) \times \Def(\Xline'),\cr
\Def(f_{\inf}) &\arr \Def(\Xline_{\inf}) \times \Def(\Xline'_{\inf}),\cr
\Def(f_{\bi}) &\arr \Def(\Xline_{\bi}) \times \Def(\Xline'_{\bi})\cr}$$
are representable by closed immersions.
\assertionitem For $(R,\alpha)$ in $\Cscr$ the inverse image of every
element in $\Def(\Xline_{\inf})(R,\alpha)$ under the canonical map
$\Def(\Xline) \arr \Def(\Xline_{\inf})$ is nonempty and a principal
homogenous space under
$$T_{\inf}(X)(R) := \Hom^{\rm Sym}_{\Zp}(T_p(X) \tensor{O_K}
T_p(X),\Gdhat_m(R)).$$
\assertionitem The inverse image of every element $\Xlinetilde_{\bi}$ in
$\Def(\Xline_{\bi})(R,\alpha)$ under the canonical map
$\Def(\Xline_{\inf}) \arr \Def(\Xline_{\bi})$ is nonempty and isomorphic to
$$T_{\bi}(X,\Xtilde_{\bi})(R) := \Hom_{O_K}(T_p(X), \Xlinetilde_{\bi}(R)).$$
\assertionitem $\Def(\Xline_{\bi})$ is trivial, i.e.\ equal to
$\Spf(W(k))$.
\assertionitem Set ${}^tf = \lambda^{-1} \circ f\vdual \circ \lambda'\colon
X' \arr X$. The inverse image of an element in $\Def(f_{\inf})(R,\alpha)$
under the canonical map $\Def(f) \arr \Def(f_{\inf})$ is non-empty and
naturally a torsor under
$$\set{(\beta,\beta') \in T_{\inf}(X)(R) \times
T_{\inf}(X')(R)}{$\beta(T_p({}^tf)(x'), x) = \beta'(x', T_p(f)(x))$}.$$
\assertionitem The inverse image of $\ftilde_{\bi} \in
\Def(f_{\bi})(R,\alpha)$ under the map $\Def(f_{\inf}) \arr
\Def(f_{\bi})$ is canonically isomorphic to
$$\set{(\gamma, \gamma') \in T_{\bi}(X,\Xtilde_{\bi})(R) \times
T_{\bi}(X,\Xtilde_{\bi})(R)}{$\ftilde_{\bi}(R) \circ \gamma = \gamma' \circ
T_p(f)$}.$$
\assertionitem $\Def(f_{\bi})$ is trivial, i.e.\ equal to $\Spf(W(k))$.}

\proof: The assertions (1), (2), (3), (5) and (6) are standard
(cf. [CN]) and hold in much greater generality. We omit the proofs.

As $\Def(\Xline_{\bi})$ is
pro-representable and formally smooth, in order to prove (4), it
suffices to show that its tangent space $t_{\Def(\Xline_{\bi})}$ is
zero-dimensional. This tangent space is canonically isomorphic to
$$\set{u \in \Hom_{O_K \tensor{\Z_p} k}(t_{X\vdual}^*,
t_{X})}{$(t_{\lambda}^{-1}) \circ (-u\star) \circ (t_{\lambda}^*) =
u$}$$
where $t_{X\vdual}^*$ denotes the $k$-dual of $t_{X\vdual}$ which is
endowed with an $O_K$-action via\break $(a\cdot t^*)(t) = t^*(\sigma(a)t)$
for $a \in O_K$, $t^* \in t_{X\vdual}^*$ and $t \in t_{X\vdual}$. The $O_K
\otimes k$-action defines decompositions $t_X = t_0 \oplus t_1$ and
$t_{X\vdual} = t'_0 \oplus t'_1$. As the Dieudonn\'e modules of
$\Xline_{\bi}$ and $\Xline_{\bi}\vdual$ are isomorphic to $\Sline^{n-2}$ we
see that $t_1 = t'_1 = 0$. By the definition of the $O_K$-action on
$t_{X\vdual}^*$ an $O_K \otimes k$-linear map $t_{X\vdual}^* \arr t_X$
corresponds to a pair of maps $(t'_1)^* \arr t_0$ and $(t'_0)^* \arr t_1$
which shows
$$\Hom_{O_K \tensor{\Z_p} k}(t_{X\vdual}^*, t_{X}) = 0.$$

Let us prove (7). As $\Xline_{\bi}$ and $\Xline'_{\bi}$ are isomorphic we
can assume that $\Xline_{\bi} = \Xline'_{\bi}$. For every $(R,\alpha)$ in
$\Cscr$ denote by $(\Xline_{\bi})_R$ the unique deformation of
$\Xline_{\bi}$ to $R$. By (1) we know that there exists at most one
deformation of $f_{\bi}$ to a morphism $(\Xline_{\bi})_R \arr
(\Xline_{\bi})_R$. Hence we have to show the existence of such a
morphism. The $p$-divisible group $\Xline_{\bi}$ is isomorphic to the
$(n-2)$-fold product of $p$-divisible $\Iscr$-modules whose Dieudonn\'e
modules are isomorpic to $\Sline$. Therefore we can assume $n = 3$. We will
now use Zink's display theory [Zi2]. The display of $(\Xline_{\bi})_R$ is
given by $(P,Q,F,V^{-1})$ where $P$ is a the free $W(R)$-module with base
$(e,f)$, where $Q$ is the submodule $\tau(1)W(R)e \oplus W(R)f$ and where
$F$ and $V^{-1}$ are given by
$$Fe = -f, \quad Ff = \tau(1)e, \quad V^{-1}\tau(1)e = -f, \quad V^{-1}f =
e.$$
The $O_K$-action is given by the decomposition $P = P_0 \oplus P_1$ with
$P_0 = W(R)e$ and $P_1 = W(R)f$, and the polarization is defined by a
perfect alternating form $\lrangle$ such that $\langle e,f \rangle =
1$. Indeed, this is clear for $R = k$. For arbitrary $R$ this is a deformation
of the display over $k$. As this deformation is unique, the given display
with $O_K$-action and polarization is the display of $(\Xline_{\bi})_R$.
To give a $p$-isogeny $(\Xline_{\bi})_R \arr (\Xline_{\bi})_R$ of
multiplicator $c$ is the same as to give a homogeneous morphism
$\varphi\colon (P,Q,F,V^{-1}) \arr (P,Q,F,V^{-1})$ such that $\langle
\varphi(x) , \varphi(y) \rangle = p^c \langle x,y \rangle$. Every such
morphism is determined by $\varphi(e)$ which has to be of the form $we$ for
some $w \in W(R)$ because of the homogenity of $\varphi$. The condition
that $\varphi$ is compatible with $F$ and $V^{-1}$ is equivalent to the
condition $w = \sigma^2(w)$ and the compatability with the form up to $p^c$
means $w\sigma(w) = p^c$. Hence it suffices to show that for any $w_0 \in
W(k)$ with these properties there exists a $w \in W(R)$ with these
properties that lifts $w_0$. As for every local Artinian ring $R$ with
residue field $k$ there is a unique homomorphism $W(k) \arr R$ we can
assume that $R = W(k)$. But then $\Delta(w_0)$ does the job, where
$\Delta\colon W(k) \arr W(W(k))$ is the homomorphism of
Cartier-Dieudonn\'e, i.e.\ the unique homomorphism of rings which is
compatible with Frobenius (see e.g.\ [Laz] VII, \pz 4).

\secstart{}\label{\uniquebiinf} Let $k$ be an algebraically closed
field of characteristic $p$ and let $R$ be a complete local noetherian
$E$-algebra with residue field $k$ and with maximal ideal $\mfr$. Let
$f_1$ and $f_2$ two $R$-valued points of $\pisog$ such that $f_1 \otimes k$
and $f_2 \otimes k$ are $\mu$-ordinary $p$-isogenies. Then the
bi-infinitesimal parts of $f_1 \otimes k$ and $f_2 \otimes k$ (in the sense
of \expodef) are isomorphic. Hence it follows from \describedef(7) that the
bi-infinitesimal parts of $f_1 \otimes R/\mfr^{n+1}$ and $f_2 \otimes
R/\mfr^{n+1}$ are isomorphic for all $n \geq 0$.

\secstart{}\label{\canonicallift} It follows that every point of
$\pisog^{\mu}$ can be lifted to characteristic zero. More precisely there
exists a canonical lift:

Let $k$ be a perfect field and let
$f\colon (A,\iota,\lambda,\hgbar) \arr (A',\iota',\lambda',\hgbar')$ be a
$\mu$-ordinary p-isogeny over $k$. It induces a homomorphism of
$p$-divisible groups
$$f[p^{\infty}]\colon X = A[p^{\infty}] \arr X' =
A[p^{\infty}].$$
This homomorphism admits a decomposition $f[p^{\infty}] =
f_{\rm et} \times f_{\rm bi} \times f_{\rm mult}$ given by the
decomposition of $X$ and $X'$ into the product of its \'etale,
bi-infinitesimal and multiplicative part. Each of the factors admits a
unique lift to an isogeny over $W(k)$: this is obvious for $f_{\rm et}$
and $f_{\rm mult}$ and it holds for $f_{\rm bi}$ because of \describedef
(7). The product of these lifts gives via the Serre-Tate theorem also a
lifting of $f$ to a $p$-isogeny $\ftilde$ over $W(k)$. This lift will be
called the {\it canonical lift} of $f$.

Its formation commutes with base change in the following
sense: If $k_1$ is a perfect extension of $k$, denote by $f_1$ the
extension of $f$ to $k_1$ and by $\ftilde_1$ the canonical lift of $f_1$ to
$W(k_1)$. Then there is a canonical identification of $\ftilde
\otimes_{W(k)} W(k_1)$ with $\ftilde_1$ where $W(k) \arr W(k_1)$ is induced
from the inclusion $k \air k_1$.

\secstart{Theorem}:\label{\muorddense} {\sl Assume that $n$ is even. Then
the $\mu$-ordinary locus of $\pisog^c \otimes \kappa(O_{E_v})$ is dense.}

\proof:  Denote by $\Mscr^{\rm nss}$ the non-supersingular locus of the
special fibre $\Mscr$ of $\Mscr_{\Iscr,K^p}$. By \sourcefinite\ we know that
the source and target morphism
$$s,t\colon Y^{\rm nss} \arr \Mscr^{\rm nss}$$
are finite. By [We1] $\Mscr^{\mu}$ is dense in $\Mscr$. Hence every component of
$Y^{\rm nss} \setminus Y^{\mu}$ has dimension strictly less than the dimension
of $\Mscr$ which is equal to $n-1$. But every component of $Y^{\rm nss}$
has dimension $\geq n-1$ \dimprop, therefore $Y^{\mu}$ has to be dense in
$Y^{\rm nss}$. It remains to show that no irreducible component of $Y$ is
contained in $Y \setminus Y^{\rm nss}$. If that were the case we would find
a point $y \in Y$ which has a neighborhood $U$ which is open in $Y$ and
which is contained in the supersingular locus $Y \setminus Y^{\rm
nss}$. Let $\Mscr^{\rm ss}$ be the supersingular locus in $\Mscr$. The
closed immersion \describedef
$$(s,t)\colon \Spf(\Oscrhat_{Y,y}) = \Spf(\Oscrhat_{U,y}) \arr
\Spf(\Oscrhat_{\Mscr,s(y)}) \times \Spf(\Oscrhat_{\Mscr,t(y)})$$
would factorize through $\Sscr = \Spf(\Oscrhat_{\Mscr^{\rm ss},s(y)}) \times
\Spf(\Oscrhat_{\Mscr^{\rm ss},t(y)})$. But by \dimsupersing\ we have
$\dim(\Sscr) = 2[(n-1)/2]$ which is strictly smaller than $n-1$ as $n$ is
even. This is a contradiction to \dimprop.

\secstart{Corollary}:\label{\topflat} {\sl Assume that $n$ is even.
\assertionlist
\assertionitem Every point of
$\pisog \otimes \kappa(E_v)$ can be lifted to characteristic zero.
\assertionitem The scheme $\pisog \otimes \kappa(E_v)$ is equi-dimensional
of dimension $n-1$.}

\proof: This follows from \muorddense, \canonicallift\ and \sourcefinite.


\paragraph{Hecke stratifications of the space of p-isogenies}

\secstart{}\label{\defineinv} Let $L$ be an unramified extension of $\Qp$
and let $K_L$ be the unique hyperspecial subgroup of $G(L)$ such that $K_L
\cap G(\Qp) = K_p$ where $K_p$ is the stabilzer of $\Lambda$ in
$V_{\Qp}$. Let $M$ and $M'$ be two $O_B$-invariant lattices of $V
\tensor{\Qp} L$ which are selfdual up to a scalar. By a lemma of Kottwitz
[Ko] there exists a $g \in G(L)$ such that $gM = M'$. The residue class of
$g$ in $K_L\backslash G(L)/K_L$ is uniquely determined and we denote it by
$\inv_L(M,M')$. For $h \in G(L)$ we have $\inv_L(hM,hM') = \inv_L(M,M')$.

Sending an $L$-rational one-parameter subgroup $\lambda$ to the double
coset of $\lambda(p)$ in $K_L\backslash G(L)/K_L$ defines a bijection between
$G(L)$-conjugacy classes of one-parameter subgroups $\G_{m,L} \arr G_L$
and $K_L\backslash G(L)/K_L$. Hence we get an identification
$$K_L\backslash G(L)/K_L = X_*(S)/\Omega(L) =
(X_*/\Omega)^{\Gal(L'/L)}$$
where $S$ is an $L$-rational maximal split torus of $G_L$ and
where $L'$ is any extension of $L$ such that $G_{L'}$ is split. Via
this identification we consider $\inv_L(M,M')$ as the unique dominant
coweight in $X_*$ which represents the corresponding element in
$X_*/\Omega$.

We endow $X_*/\Omega$ with a partial order by saying $\lambda\Omega \leq
\lambda'\Omega$ if $\lambda\Omega$ is in the convex hull of
$\lambda'\Omega$ in $X_* \otimes \R$.

\secstart{} Let $\sbar \arr \pisog$ be a geometric point in characteristic
zero corresponding to a $p$-isogeny $f\colon (A,\iota,\lambda,\hgbar) \arr
(A',\iota',\lambda',\hgbar')$. We get an induced map $T_p(f)\colon T_p(A)
\arriso T_p(A')$. Due to the determinant condition we can choose an
$O_B$-linear symplectic similitude $T_p(A) \arriso \Lambda$ (e.g.\ [RZ]
6.10) and via this isomorphism we consider $T_p(A)$ and $T_p(A')$ as
lattices in $V_{\Qp}$. Both are $O_B$-invariant and selfdual up to a scalar.
Let $s \in S$ be the topological image of $\sbar$. Then
$$\eta^0(s) := \inv_{\Qp}(T_p(A),T_p(A'))$$
is independent of the choice of $\sbar$ over $s$ and the choice of the
isomorphism $T_p(A) \cong \Lambda$. Altogether we obtain a map
$$\eta^0\colon \pisog \otimes E \arr (X_*/\Omega)^{\sigma}$$
where $\sigma$ denotes the Frobenius automorphism. This map is clearly locally
constant. Hence $\pisog \otimes E$ is the union of open and closed subschemes
$\pisog^{\alpha}_E$ with\break $\alpha \in (X_*/\Omega)^{\sigma}$.

\secstart{}\label{\stratmuord} Let $\pisog^{\mu}$ be the $\mu$-ordinary
locus of $\pisog \otimes \kappa(E_v)$ and $s$ be a point of
$\pisog^{\mu}$. Let $\sbar\colon \Spec(k) \arr \pisog^{\mu}$ be a geometric
point over $s$ corresponding to a $p$-isogeny $f\colon
(A,\iota,\lambda,\hgbar) \arr (A',\iota',\lambda',\hgbar')$ over $k =
\kappa(\sbar)$.

By \canonicallift\ there exists a canonical lift $\ftilde$ of $f$ to
$W(k)$. Denote by
$$\alpha = T_p(\ftilde^0)\colon (T_p(A^0),\iota,\lrangle) \arr
(T_p({A'}^0),\iota',\lrangle')$$
the map on Tate modules which is induced by the generic fibre of
$\ftilde$. Because of the determinant condition we can
choose a $B$-linear symplectic similitude of\break $(T_p(A^0),\iota,\lrangle)
\otimes \Qp  = (T_p({A'}^0),\iota',\lrangle') \otimes \Qp)$ with the
skew-hermitian $B$-module $(V_{\Qp}, \lrangle_{\Qp})$. Hence we can
consider $T_p(A^0)$ and $T_p({A'}^0)$ as $O_B$-lattices in $V_{\Qp}$ which
are selfdual up to a scalar.

By construction the induced homomorphism of $\ftilde$ on $p$-divisible groups
splits into a product of a lift of the \'etale, the bi-infinitesimal, and
the multiplicative part of $f$. Hence there exists an $h \in L(\Qp)$ such
that $hM = M'$ where $L \subset G_{\Qp}$ is the Levi subgroup defined in
\defineLevi. Denote by ${}_LK$ the intersection $K \cap L(\Qp)$. The double
coset
$$\eta^{\mu}(s) \in {}_LK\backslash L(\Qp)/{}_LK$$
of $h$ is well defined and depends only on $s$. As in \defineinv\ we can
consider $\eta^{\mu}(s)$ as an element of ${}_LX_*/{}_L\Omega$ where
$({}_LX^*,{}_LR^*,{}_LX_*,{}_LR_*,{}_L\Delta)$ is the based root datum of
$L$ with Weyl group ${}_L\Omega$. In fact $\eta^{\mu}(s)$ lies
automatically in the $\sigma$-invariants of ${}_LX_*/{}_L\Omega$. Note
further that $\eta^{\mu}(s)$ does not depend on the choice of $\sbar$
\canonicallift.

\secstart{Proposition}:\label{\decompmuord} {\sl The map
$$\eta^{\mu}\colon \pisog^{\mu} \arr {}_LX_*/{}_L\Omega$$
is locally constant.}

\proof: It suffices to show the following: Let $g\colon
(A,\iota,\lambda,\hgbar) \arr (A',\iota',\lambda',\hgbar')$ be an
$S$-valued point of $\pisog^{\mu}$ where $S$ is a $\kappa(O_{E_v})$-scheme
which is of the form $S = \Spec(k\dlbrack t \drbrack)$ for some
algebraically closed field $k$. Let $s$ be the special point of $S$ and $t$
be a geometric pont of $S$ lying over the generic point of $S$. Then
$\eta^{\mu}(s) = \eta^{\mu}(t)$. Denote by $f\colon \Xline \arr \Xline'$
the associated morphism of $p$-divisible $\Iscr$-modules. Denote by $G'$
the reductive $\Qp$-group $G' = GL_{O_K}(V_{\Qp})$ and by $L'$ the
centralizer of $N(\mu)$ \defineLevi\ in $G'$. The inclusion $L \arr
L'$ induces a map ${}_LX_*/{}_L\Omega \arr {}_{L'}X_*/{}_{L'}\Omega$. It is
easy to see that this map is injective. Hence it suffices to prove that the
type of $f$ as a $p$-isogeny of $p$-divisible $O_K$-modules is locally
constant. As the \'etale and the multiplicative rank of $X$ and $X'$
are locally constant functions on $S$ we get a filtration of
$p$-divisible $O_K$-modules
$$0 \subset (X,\iota)_{\mult} \subset (X,\iota)_{\inf} \subset (X,\iota)$$
whose successive quotients are $(X,\iota)_{\bi}$ resp.\ $(X,\iota)_{\et}$,
ditto for $(X',\iota')$ ([Gr] III, 7.4). We obtain an induced $p$-isogeny on
the associated graded $p$-divisible $O_K$-modules
$${\rm gr}(f) = f_{\mult} \times f_{\bi} \times f_{\et}$$
and by definition the type of $f$ and that of ${\rm gr}(f)$ define the same
function on $S$. Hence we can assume that $X$ and $X'$ are products of
their multiplicative, their bi-infinitesimal and their \'etale part.

We set $k_1 = \kappa(t)$. Let $\ftilde_s$ (resp.\ $\ftilde_t$) be the
canonical lift of $f_t$ (resp.\ $f_s$) to $W(\kappa(s)) = W(k)$ (resp.\ to
$W(k_1)$). It suffices to find a local ring $R$ with $R/pR = k\dlbrack t
\drbrack$ together with $W(k)$-homomorphisms
$R \arr W(k)$ and $R \arr W(k_1)$ and a $O_K$-linear $p$-isogeny $\ftilde$
over $R$ which lifts $f$ such that $\ftilde \otimes W(k) = \ftilde_s$ and
$\ftilde \otimes W(k_1) = \ftilde_t$. For $R$ we take $W(k)\dlbrack T
\drbrack$. As homomorphism $R \arr W(k)$ (resp. $R \arr W(k_1)$) we take
the unique $W(k)$-homomorphism which sends $T$ to zero (resp.\ to the
Teichm\"uller representative of $t$ in $W(k_1)$). Our $\ftilde$ will be a
product $\ftilde = \ftilde_{\et} \times \ftilde_{\bi} \times
\ftilde_{\mult}$ where $\ftilde_{?}$ lifts $f_{?}$. To lift $f_{\et}$ and
$f_{\mult}$ is trivial, hence we have only to construct
$\ftilde_{\bi}$. Because of the uniqueness of lifts of $f_{s,\bi}$ and
$f_{t,\bi}$ \describedef, it suffices to construct a compatible system of
lifts $\ftilde_{\bi,n}$ of $f_{\bi}$ to $W_n(k)\dlbrack t \drbrack$. This
can be done as in the proof of \describedef(7) using the theory of
Dieudonn\'e windows developed in [Zi3] (one can take as a frame $(A =
W(k)\dlbrack t \drbrack, p^nA, \sigma)$ where $\sigma$ is the unique
continuous ring endomorphism of $A$ with $\sigma(T) = T^p$ which induces on
$W(k)$ the Frobenius).

\secstart{}\label{\Tatereduction} Let $k$ be a perfect field of
characteristic $p$ and let $f\colon (X,\iota,\lambda) \arr
(X',\iota',\lambda')$ be a homomorphism of $p$-divisible $\Iscr$-modules
which is induced by a $\mu$-ordinary $p$-isogeny over $k$. Denote by
$\ftilde^0$ the generic fibre of the canonical lift of $f$ to $W(k)$. It is
the product $\ftilde^0_{\et} \times \ftilde^0_{\bi} \times
\ftilde^0_{\mult}$ where $\ftilde^0_{?}$ is the generic fibre of the lift
of $f_{?}$. Then the invariant factors of $T_p(\ftilde^0_{\et})$ (resp.\ of
$T_p(\ftilde^0_{\mult})$) in the sense of [BouA] chap.\ VII, \pz 4, n${}^o$
6 are the same as the invariant factors of $T_p(f)$ (resp.\ of
$T_p({}^tf)$). This follows from (duality and):

\claim Lemma: Let $S$ be a connected noetherian scheme and let $f\colon X
\arr X'$ be an isogeny of $p$-divisible groups. Assume that $X_s[p]$
is an \'etale group scheme for all $s \in S$. Then the invariant factors of
$T_p(f_{\sbar})$ are the same for all geometric points $\sbar \arr S$.

\proof: The functor $X \asr T_p(X_{\sbar})$ sets up an equivalence between the
category of \'etale $p$-divisible groups over $S$ 
and the category of finitely generated free $\Zp$-modules with continuous
$\pi_1(S,\sbar)$-action and the latter category is independent of the
choice of $\sbar$ up to inner automorphism.

\secstart{} Altogether we obtain two decompositions
$$\eqalign{\pisog \otimes E &= \coprod_{\alpha \in (X_*/\Omega)^{\sigma}}
\pisog^{\alpha}_E, \cr
\pisog^{\mu} &= \coprod_{\gamma \in ({}_LX_*/{}_L\Omega)^{\sigma}}
\pisog^{\mu,\gamma} \cr}$$
into open and closed subschemes.

\secstart{Proposition}:\label{\sourcelocallyfree} {\sl Source and target
morphism $\pisog^{\mu} \arr \Mscr^{\mu}$ are finite locally free.}

\proof: It suffices to prove this for the source morphism $s$ (by duality). We
know already that the source morphism is finite \sourcefinite. For every
point $x$ of $\Mscr$ the fibre $s^{-1}(\xbar)$ of a geometric point $\xbar$
lying over $x$ is the spectrum of a finite $\kappa(\xbar)$-algebra
$A_{\xbar}$. As $\Mscr^{\mu}$ is regular and in particular reduced, it
suffices to show that the function $x \asr \dim_{\kappa(\xbar)}(A_{\xbar})$
is locally constant. From \finitefibres\ we already know that $x \asr
\#s^{-1}(\xbar)$ is locally constant. Because of \decompmuord\ it suffices
to show: Let $\ybar \arr \pisog^{\mu,\gamma}$ be a geometric
point and $\xbar = s(\ybar)$ its image in $\Mscr^{\mu}$. Then the
$\kappa(\xbar)$-dimension of $\Oscr_{s^{-1}(\xbar),\ybar} =
\hat{\Oscr}_{s^{-1}(\xbar),\ybar}$ depends only on $\gamma$. But it follows
from \describedef\ and \uniquebiinf\ that the isomorphism class of the
$\kappa(\xbar)$-algebra $\hat{\Oscr}_{s^{-1}(\xbar),\ybar}$ depends only on
$T_p(f)$ and $T_p({}^tf)$ and in particular only on $\gamma$ \Tatereduction.

\secstart{}\label{\definereduction} For any locally noetherian scheme $Y$
we denote by $\Q[Y]$ the $\Q$-vector space which has as a basis the
components of $Y$. If $V \subset Y$ is an open dense subset we have a
canonical identification $\Q[V] = \Q[Y]$ of $\Q$-vector spaces.

We will consider the cases $Y = \pisog \otimes E$ and $Y = \pisog \otimes
\kappa(E_{\nu})$. Reduction of cycles defines a map
$${\rm red}_p\colon \Q[\pisog \otimes E] \arr \Q[\pisog \otimes
\kappa(E_v)]$$
which associates to each component of $\pisog \otimes E$ the intersection
of its closure in $\pisog$ with $\pisog \otimes \kappa(E_v)$

\secstart{}\label{\definecomposition} Set $\pisog^0 = \pisog \otimes E$ and
$\Mscr^0 = \Mscr \otimes E$. Consider the morphism
$${\rm comp}\colon \pisog^0 \times_{t,\Mscr^0,s} \pisog^0 \arr \pisog^0,
\quad (f,g) \asr g \circ f$$
which is given by composition of $p$-isogenies. This morphism is finite
\'etale because source and target morphism $s,t\colon \pisog^0 \arr
\Mscr^0$ are finite \'etale, t. If $Z_1$ and $Z_2$ are two
irreducible components of $\pisog^0$, we form the fibre product of
$$\matrix{&&\Mscr^0\cr
&&\addright{\Delta_{\Mscr^0}}\cr
Z_1 \times Z_2 & \airrover{t,s} & \Mscr^0 \times \Mscr^0.\cr}$$
This is a cycle in codimension zero of $\pisog^0 \times_{t,\Mscr^0,s}
\pisog^0$ and its push forward via ${\rm comp}$ is an element in
$\Q[\pisog^0]$. This defines the structure of a $\Q$-algebra on
$\Q[\pisog^0]$.

In the same way (by replacing $?^0$ by $?^{\mu}$ and using
\sourcelocallyfree) we get a multiplication on $\Q[\pisog^{\mu}]$.
If $n$ is even, \muorddense\ implies that we have $\Q[\pisog \otimes
\kappa(E_{\nu})] = \Q[\pisog^{\mu}]$. Hence $\Q[\pisog \otimes
\kappa(E_{\nu})]$ obtains the structure of a $\Q$-algebra such that ${\rm
red}_p$ \definereduction\ becomes a homomorphism of $\Q$-algebras.

\secstart{}\label{\maintheorem} Assume that $n$ is even. In [We2] \pz 4 a
$\Q$-algebra homomorphism
$$h^0\colon \Hscr_{\Q}(G(\Qp)//K_p)_0 \arr \Q[\pisog \otimes E]$$
is constructed which sends the characteristic function of a $K_p$-double
coset $\alpha$ of $G(\Q_p)$ to the corresponding stratum
$(\pisog \otimes E)^{\alpha}$. This stratum consists of those points $f$ in
$\pisog \otimes E$ where $f$ is given on Tate modules by $\alpha$ using the
identification \Cartanident\ (see loc.\ cit.\ for the precise definition).

In $\Q[\pisog \otimes \kappa(E_v)]$ we have the element ${\rm Frob}$
corresponding to relative Frobenius isogenies as in \defineFrob.

The principal result is now:

\claim Main Theorem: We consider $\Q[\pisog \otimes
\kappa(E_v)]$ via ${\rm red}_p \circ h^0$ as an algebra over
$\Hscr_{\Q}(G)$. Then we have
$$H_{G,\{\mu\}}({\rm Frob}) = 0.$$

\secstart{}\label{\defineHeckealg2} To prove the main theorem we use the
stratification of $\pisog^{\mu}$ indexed by elements in
${}_LX_*/{}_L\Omega$ and \stratmuord. We define a $\Q$-algebra homomorphism
$$\hbar\colon \Hscr_{\Q}(L)_0 \arr \Q[\pisog^{\mu}] = \Q[\pisog \otimes
\kappa(E_v)]$$
as in [We2] \pz 5 where $\Hscr_{\Q}(L)_0$ denotes again the subalgebra of
functions with support in $\End(\Lambda)$. It is defined as the map which
sends the characteristic function of an ${}_LK$-double coset (corresponding
to a $\gamma \in {}_LX_*/{}_L\Omega$) to the cycle of the
corresponding stratum divided by a constant $a(\gamma)$. This constant is
defined as the formal multiplicity of the geometric fibres of the source
morphism $\pisog^{\mu} \arr \Mscr^{\mu}$. It depends only on
$\gamma$ by \sourcelocallyfree.

Further denote by
$$\Sscrdot\colon \Hscr_{\Q}(G) \arr \Hscr_{\Q}(L)$$
the twisted Satake homomorphism defined in [We2] \pz 1. As in [FC] VII,4 we
have a commutative diagram
$$\matrix{\Hscr_{\Q}(G)_0 & \arrover{h^0} & \Q[\pisog \otimes E] \cr
\addleft{\Sscrdot} && \addright{{\rm red}_p} \cr
\Hscr_{\Q}(L)_0 & \arrover{\hbar} & \Q[\pisog \otimes
\kappa(E_v)].\cr}\lformno\sublabel{\commdiag}$$

\secstart{Proof of the main theorem \maintheorem}: Recall that we have
defined in \defineLevi\ an element $\{N(\mu)\} \in
(X_*/{}_L\Omega)^{\sigma}$. This is the same as the element $[\mu]$ defined
in [We2] 2.4. Let $g \in \Hscr_{\Q}(L)$ be the corresponding characteristic
function \defineHeckealg2. Let $H_{G,\{\mu\}} \in \Hscr_{\Q}(G)[X]$ be the
Hecke polynomial defined in loc.\ cit.\ 2.6. We regard $\Hscr_{\Q}(L)$ via
the twisted Sataka homomorphism $\Sscrdot$ as an algebra over
$\Hscr_{\Q}(G)$. By loc.\ cit.\ 2.9 we have
$$H_{G,\{\mu\}}(g) = 0.$$
Hence by \commdiag\ it suffices to show that the element $g$ is sent to
the element ${\rm Frob}$ in $\Q[Y \otimes \kappa(E_v)]$ via $\hbar$ where
${\rm Frob}$ is the relative Frobenius isogeny with respect to $\kappa(E_v)$.

For this we use the notations of \moreno. We may assume (to simplify
notations) that $n > 2$, hence $n
\geq 4$ as $n$ is even, as the case $n = 2$ is the classical case (of
elliptic curves). Then $E$ is a quadratic imaginary extension of $\Q$ and
we have $E_v = K$. We have $B \otimes \Q_p = M_m(K)$ und using Morita
equivalence we can assume that $m = 1$. Hence we can identify $G_{\Qp}$
with $GU(V_{\Qp}, \beta)$ where $\beta\colon V_{\Qp} \times V_{\Qp} \arr K$
is the skew-hermitian form such that its composition with $\Tr_{K/\Qp}$ is
$\lrangle_{\Qp}$. Using the canonical embedding $G_{\Q_p} \air GSp(V_{\Qp},
\lrangle_{\Qp})$ we identify $X_*$ with
$$\set{(x_i) \in \Z^{2n}}{$x_i + x_{2n+1-i} = {\rm const}$ for all $i$}$$
and $\Omega$ with $S_n$ which acts on $X_*$ by permuting simlutaneously
$(x_1,\ldots,x_n)$ and $(x_{2n},\ldots,x_{n+1})$. The element $\{\mu\} \in
X_*/\Omega$ is given by the class of $(1,0^{n-1},1^{n-1},0)$ because of our
assumption on the signature. Further the nontrivial element $\sigma$ in
$\Gal(K/\Q_p)$ acts by sending $(x_i)_{1\leq i \leq 2n}$ to $(x_{n+1},
\ldots, x_{2n}, x_1, \ldots, x_n)$. Hence the norm of $\{\mu\}$ is given by
$$N(\{\mu\}) = (2,1^{n-2},0,2,1^{n-2},0) \in X_*/{}_L\Omega.$$
But this is precisely the type of $V^2$ on a $\mu$-ordinary unitary
Dieudonn\'e module \muorddieudmod. As ${\rm Frob}$ corresponds to $V^2$ via
covariant Dieudonn\'e theory, we see that ${\rm Frob}$ is the
cycle in $\Q[\pisog \otimes \kappa(E_v)]$ corresponding to $\gamma =
N(\{\mu\})$. Further the definition of $a(\gamma)$ \defineHeckealg2\ shows that
$a(N(\{\mu\})$ is equal to one which proves the claim.


\bibliography

\indention{[BouA]\ }
\litem{[BR]} D. Blasius, J.D. Rogawski: {\it Zeta functions of Shimura
  varieties}, in:\break U. Jannsen, S. Kleiman, J.-P. Serre (ed.): {\it
Motives}\/ 2, Proc. Symp. Pure Math. {\bf 55} (1994), 447--524.
\litem{[BouA]} N. Bourbaki: {\it Algebra}, chap. IV-VII, Springer (1988).
\litem{[De]} P. Deligne: {\it Vari\'et\'es de Shimura}, in A. Borel,
  W. Casselman (ed.): {\it Automorphic forms, representations, and
  $L$-functions}\/ 2, Proc. Symp. Pure Math. {\bf 33} (1977), 247--289.
\litem{[Gr]} A.~Grothendieck: ``Groupes de Barsotti-Tate et cristaux de
Dieudonn\'e,'' Presses de l'Universit\'e de Montr\'eal, 1974.
\litem{[Ill]} L. Illusie: {\it Complexe de de\thinspace Rham-Witt et
cohomologie cristalline}, Ann. Sci. \'Ecole Norm. Sup. (4) {\bf 12} (1979),
501--661.
\litem{[Ko]} R. Kottwitz: {\it Points on some Shimura varieties over finite
fields}, J. AMS {\bf5} (1992), 373-444.
\litem{[Kr]} H. Kraft: {\it Kommutative algebraische $p$-Gruppen (mit
Anwendungen auf $p$-divisible Gruppen und abelsche Variet\"aten)},
lecture manuscript, Universit\"at Bonn, 1975.
\litem{[Laz]} M. Lazard: {\it Commutative formal groups},
Lecture Notes in Mathematics {\bf 443}, Springer-Verlag, Berlin-New York, 1975.
\litem{[LM]} G. Laumon, L. Moret-Bailly: ``Champs alg\'ebriques'',
Ergebnisse der Mathematik und ihrer Grenzgebiete (3. Folge) {\bf 39},
Springer (2000).
\litem{[Mo]} B. Moonen: {\it Group schemes with additional structure and Weyl
group cosets}, in ``Moduli of Abelian Varieties'' (ed. by
C. Faber, G. van der Geer, F. Oort), Progress in Mathematics {\bf 195},
Birkh\"auser (2001).
\litem{[Ra]} M. Raynaud: {\it Sch\'emas en groupes de type $(p,\dots, p)$},
Bull. Soc. Math. France {\bf 102} (1974), 241--280.
\litem{[RZ]} M. Rapoport, T. Zink: {\it Period Spaces for $p$-divisible
Groups}, Annals of Mathematics Studies {\bf141}, Princeton University Press,
Princeton (1996).
\litem{[We1]} T. Wedhorn: {\it Ordinariness in good reductions of Shimura
varieties of PEL-type}, Ann. Sci. de l'ENS {\bf 32} (1999), 575 - 618.
\litem{[We2]} T. Wedhorn: {\it Congruence relations for some Shimura
  varieties}, J. reine angewandte Math. (Crelle) {\bf 524} (2000), 43--71.
\litem{[We3]} T. Wedhorn: {\it The dimension of Oort strata of Shimura
varieties of PEL-type}, in ``Moduli of Abelian Varieties'' (ed. by
C. Faber, G. van der Geer, F. Oort), Progress in Mathematics {\bf 195},
Birkh\"auser (2001).
\litem{[Zi1]} T. Zink: {\it Cartiertheorie kommutativer formaler Gruppen},
Teubner-Texte zur Mathematik {\bf 68}, Leipzig (1984).
\litem{[Zi2]} T. Zink: {\it Cartier Theory and Crystalline Dieudonn\'e
Theory}, preprint Bielefeld (1998).
\litem{[Zi3]} T. Zink: {\it Windows for Displays of $p$-divisible Groups},
in ``Moduli of Abelian Varieties'' (ed. by C. Faber, G. van der Geer,
F. Oort), Progress in Mathematics {\bf 195}, Birkh\"auser (2001).

\bye